\documentclass[12pt]{amsart}
\usepackage[T1]{fontenc}
\usepackage{lmodern,amsmath,amssymb,mathtools,mathrsfs,booktabs,microtype}
\usepackage[margin=1in]{geometry}
\usepackage[hidelinks]{hyperref}
\newtheorem{theorem}{Theorem}[section]
\newtheorem{lemma}[theorem]{Lemma}
\newtheorem{proposition}[theorem]{Proposition}
\newtheorem{corollary}[theorem]{Corollary}
\theoremstyle{remark}
\newtheorem{remark}[theorem]{Remark}
\newtheorem{example}[theorem]{Example}
\newcommand{\PP}{\mathbb P}
\newcommand{\CC}{\mathbb C}
\newcommand{\ZZ}{\mathbb Z}
\newcommand{\QQ}{\mathbb Q}
\newcommand{\OO}{\mathcal O}
\newcommand{\HH}{\mathcal H}
\newcommand{\LL}{\mathcal L}
\newcommand{\MC}{\mathcal C}
\newcommand{\MM}{\mathcal M}
\newcommand{\YY}{\mathcal Y}
\DeclareMathOperator{\Sym}{Sym}
\DeclareMathOperator{\Gr}{Gr}
\DeclareMathOperator{\Spec}{Spec}
\DeclareMathOperator{\Ext}{Ext}
\DeclareMathOperator{\Ob}{Ob}
\DeclareMathOperator{\rank}{rank}
\DeclareMathOperator{\coker}{coker}
\DeclareMathOperator{\tr}{tr}
\DeclareMathOperator{\edim}{edim}
\DeclareMathOperator{\vdim}{vdim}
\DeclareMathOperator{\Aut}{Aut}
\newcommand{\Mbar}{\overline M}
\title[Conics and double covers]{Enumerative geometry of conics and double covers on projective complete intersections}
\author{Dang Tuan Hiep}
\address{Faculty of Mathematics and Computer Science, Da Lat University, Vietnam}
\email{hiepdt@dlu.edu.vn}
\author{Ha Thi Thao Tram}
\address{Faculty of Science, Ho Chi Minh City University of Agriculture and Forestry, Vietnam}
\email{htthaotram@hcmuaf.edu.vn}
\date{\today}
\begin{document}
\begin{abstract}
We compare genus zero degree two Gromov--Witten invariants with the number of embedded conics on general projective complete intersections whose conic space has expected dimension zero. Their difference comes from double covers of lines. We compute this contribution as an integral of Chern classes over the Fano scheme of lines, using an explicit obstruction bundle on the compactified space of double covers.
\end{abstract}
\maketitle

\section{Introduction}
Throughout, we work over $\CC$. Let
\[
X=X_{m_1,\ldots,m_s}\subset\PP^r,\quad
m=\sum_{i=1}^s m_i,\quad n=r-s,
\]
be the zero locus of general equations of positive degrees $m_i$. Here general means that the equations belong to a Zariski open subset on which $X$ is smooth and the conclusions of Lemma~\ref{lem:geometry} hold simultaneously. The degree of a stable map $g:C\to X$ is
$\deg g^*\OO_X(1)$. The degree two invariant is summed over the curve classes with hyperplane degree two if more than one such class occurs. The condition
\begin{equation}\label{eq:dimension}
2m+s = 3r-1
\end{equation}
is equivalent to the vanishing of the virtual dimension of the Kontsevich moduli space of genus zero stable maps of degree two 
\[
\vdim\Mbar_{0,0}(X,2)=n-3+2(r+1-m) = 0.
\]
It also makes the expected dimension of the Hilbert scheme of conics zero and forces
\[n = 2(m-r)+1\] 
to be odd.

Let $N_2(X)$ denote the genus zero degree two Gromov--Witten invariant, and let $n_2(X)$ denote the number of embedded conics. These numbers need not agree, since stable maps of degree two also include double covers of lines. When $n=3$, the lines are isolated, and each contributes $1/8$, giving
\[
N_2(X)=n_2(X)+\frac18n_1(X),
\]
where $n_1(X)$ is the number of lines on $X$. When $n\geq5$ and $X$ contains a line, the lines move in a family of dimension 
\[h=\frac{n-3}{2}.\] 
Their compactified double covers form a family of dimension $h+2$, whose contribution to $N_2(X)$ is computed using its obstruction theory.

Let $G_1=\Gr(2,r+1)$ be the Grassmannian parametrizing lines in $\PP^r$, with tautological exact sequence
\[
0\longrightarrow U
\longrightarrow \CC^{r+1}\otimes\OO_{G_1}
\longrightarrow Q
\longrightarrow 0,
\]
where $U$ and $Q$ have ranks $2$ and $r-1$, respectively. The defining equations of $X$ determine a section of
\[
A_X=\bigoplus_{i=1}^s\Sym^{m_i}U^\vee,
\]
whose zero scheme is the Fano scheme of lines 
\[F=F_1(X)\subset G_1.\] 
Its expected dimension is
\[
h =2r-2-(m+s)=\frac{n-3}{2}.
\]
For general defining equations, $F$ is empty or smooth of pure dimension $h$, as shown in Lemma~\ref{lem:geometry}. 

Let $f_1,\ldots,f_s$ be the defining equations of $X$, with $\deg f_i=m_i$. For a line $\ell\subset X$, differentiation in directions normal to $\ell$, followed by restriction to $\ell$, gives a linear map
\[
\delta_\ell:Q_\ell\longrightarrow\bigoplus_{i=1}^sH^0(\ell,\OO_\ell(m_i-1)).
\]
These maps assemble into a bundle homomorphism
\begin{equation}\label{eq:delta}
\delta:Q|_F
\longrightarrow
\bigoplus_{i=1}^s\Sym^{m_i-1}U^\vee|_F.
\end{equation}
Lemma~\ref{lem:normal} shows that
\[
\ker\delta_\ell
=
H^0(\ell,N_{\ell/X}(-1))=0,
\]
for every $\ell\in F$. Thus $\delta$ has constant rank
$r-1$, and its cokernel is locally free. Define
\begin{equation}\label{eq:B}
B=\det U^\vee\otimes\coker(\delta).
\end{equation}
Since 
\[\rank\Sym^{m_i-1}U^\vee=m_i,\] 
we have
\[
\rank B=m-r+1=h+2,
\]
where the last equality follows from \eqref{eq:dimension}.
Restrictions to $F$ are understood whenever the base is clear.

Let $x,y$ be the Chern roots of $U^\vee$, and put
\[
\sigma = x + y = c_1(U^\vee),
\quad
Y_1=\{2x,x+y,2y\}.
\]
Thus $Y_1$ is the multiset of Chern roots of
$\Sym^2U^\vee$.
For a vector bundle $E$, write
\[
c_t(E)=\sum_{j\geq0}c_j(E)t^j,
\quad
c_{t/2}(E)=\sum_{j\geq0}c_j(E)(t/2)^j.
\]

\begin{theorem}\label{thm:main}
Assume \eqref{eq:dimension}, $n\geq3$, and general defining equations. The Hilbert scheme of conics on $X$ is finite and reduced, and it contains only smooth conics. The Fano scheme $F=F_1(X)$ is either empty or smooth of pure dimension $h$. Then
\begin{equation}\label{eq:main}
N_2(X)=n_2(X)+D_2(X),\quad
D_2(X)=\frac18\int_F[t^h]\frac{c_t(B)}{c_{t/2}(\Sym^2U^\vee)}.
\end{equation}
Here $[t^h]$ denotes the coefficient of $t^h$ in the formal power series expansion at $t=0$. The integral is understood to be zero if $F$ is empty. In particular, when $n=3$,
\begin{equation}\label{eq:isolated}
N_2(X)=n_2(X)+\frac18\#F.
\end{equation}
\end{theorem}

The correction is obtained from an obstruction bundle defined on the entire compactified space of double covers, including the nodal boundary. The definition of $B$ gives
\begin{equation}\label{eq:Kclass}
[B]
=
\sum_{i=1}^s
[\det U^\vee\otimes\Sym^{m_i-1}U^\vee]
-
[\det U^\vee\otimes Q]
\qquad\text{in }K^0(F).
\end{equation}
Hence $c_t(B)$ is the restriction to $F$ of the class
\[
\frac{\displaystyle\prod_{i=1}^s
c_t(\det U^\vee\otimes\Sym^{m_i-1}U^\vee)}
{c_t(\det U^\vee\otimes Q)}
\]
on $G_1$. Its expression therefore depends only on the multidegree of $X$, rather than on the chosen defining equations. 

Under the hypotheses of Theorem~\ref{thm:main}, we have
\begin{equation}\label{eq:ambient}
D_2(X)
=
\frac18\int_{G_1}[t^h]
\frac{\displaystyle\prod_{i=1}^s
c_t(\det U^\vee\otimes\Sym^{m_i-1}U^\vee)}
{c_t(\det U^\vee\otimes Q)
 c_{t/2}(\Sym^2U^\vee)}c_{m+s}(A_X).
\end{equation}
Define the polynomial of degree $h$
\begin{equation}\label{eq:Kpoly}
K_X(x,y)
=
[t^h]\,
\frac{\displaystyle
\prod_{i=1}^s\prod_{j=0}^{m_i-1}
\bigl(1+((j+1)x+(m_i-j)y)t\bigr)}
{(1+(x+y)t)^{r+1}(1+(x+y)t/2)}.
\end{equation}
Evaluating $K_X(x,y)$ at the Chern roots of $U^\vee$ gives
the following equality 
\[
K_X(x,y)
=
[t^h]\,
\frac{\displaystyle\prod_{i=1}^s
c_t(\det U^\vee\otimes\Sym^{m_i-1}U^\vee)}
{c_t(\det U^\vee\otimes Q)
 c_{t/2}(\Sym^2U^\vee)} \in A^h(G_1)_{\mathbb Q}.
\]
Thus, $K_X(x,y)$ is an explicit polynomial representative
of the degree $h$ Chow class occurring in
\eqref{eq:ambient}.
\begin{theorem}\label{thm:coefficient}
With the above notation, we have
\begin{equation}\label{eq:ranktwo}
D_2(X)
=
\frac18[x^ry^{r-1}]
\left(
(x-y)K_X(x,y)
\prod_{i=1}^s\prod_{j=0}^{m_i}
(jx+(m_i-j)y)
\right).
\end{equation}
Here $[x^ry^{r-1}]$ denotes the coefficient of $x^ry^{r-1}$ in the polynomial inside the parentheses.
\end{theorem}

For isolated lines, the correction is given by the Aspinwall--Morrison multiple cover formula~\cite{AM}, proved mathematically by Voisin~\cite{Voisin}. Faber--Pandharipande~\cite[Section 4.1]{FP} computed the corresponding local integrals by virtual localization. A degree $d$ cover contributes $1/d^3$, giving $1/8$ in degree two. Our formula recovers this case and expresses the contribution of double covers of moving lines in terms of Chern classes on the Fano scheme of lines.

The invariant $N_2(X)$ can also be computed using Kontsevich's localization formula and virtual localization~\cite{Kontsevich,GP}. For complete intersection Calabi--Yau threefolds, the graph sums were evaluated through degree six in~\cite{Hiep2016}. The five degree two values in Section~\ref{sec:examples} reproduce these known results and serve as checks of our formula. Katz~\cite[Section 3]{Katz} constructed the projective bundle of plane quadrics and the rank eleven quotient bundle used to count conics on a quintic threefold. 

Bonavero--H\"oring~\cite[Theorem 1.1]{BH} counted conics through two general points on general complete intersections satisfying $2m=r+s+1$. Their conic space has the expected dimension $2n-2$, and a general conic is a quasi-line. Their Section~2.3 also provides the quotient bundle construction and the global generation argument recalled below. Our condition $2m+s=3r-1$ instead makes the unconstrained conic space zero-dimensional. In this setting, we compare its enumerative count with the degree two Gromov--Witten invariant by computing the contribution of double covers of lines.

Ciliberto--Zaidenberg~\cite{CZ2020} studied enumerative formulas and numerical invariants of Fano schemes of complete intersections. Their survey~\cite[Sections 2--4]{CZ2022} provides a broader account of lines and conics. Our coefficient formulas follow from Gysin formulas for Grassmannians and flag bundles~\cite{DP,Hiep2019,HGV}. For background on intersection theory and enumerative geometry, we refer to Fulton~\cite{FultonIT} and Eisenbud--Harris~\cite{EH3264}. Schubert classes and degeneracy loci are treated in Fulton--Pragacz~\cite{FultonPragacz}.

Jinzenji--Nakamura--Suzuki~\cite[Theorem 5.1]{JNS} computed the contribution of double covers of moving lines on projective hypersurfaces. Their construction includes nodal source curves and a square root line bundle on the compactified cover stack. In the hypersurface case, our correction formula is equivalent to theirs. We extend this formula to complete intersections under \eqref{eq:dimension} and provide a description of families, along with a scheme-theoretic separation of the double cover component.

Moving multiple covers have also been studied in the enumerative geometry of higher-dimensional Calabi--Yau varieties~\cite{PZ}. In our setting,
\[
r+1-m=\frac{3-n}{2}.
\]
Thus, by adjunction,
\[
\omega_X\simeq\OO_X\left(\frac{n-3}{2}\right).
\]
In particular, $X$ has ample canonical bundle and is
of general type when $n>3$.

Our main result is the correction formula \eqref{eq:main} for double covers of moving lines and its expression \eqref{eq:ambient} as an integral over $G_1$. We construct the obstruction bundle on the entire compactified cover space, including the nodal boundary, and prove the branch gerbe description for families over arbitrary bases. The stable map localization formula and the Hilbert scheme formula provide independent numerical checks.

\subsection*{Outline}
Section~\ref{sec:geometry} establishes the geometric properties of lines and conics. Section~\ref{sec:covers} constructs the obstruction bundle and proves Theorems~\ref{thm:main} and~\ref{thm:coefficient}. Sections~\ref{sec:localization} and~\ref{sec:hilbert} give formulas for $N_2(X)$ and $n_2(X)$, respectively. Section~\ref{sec:examples} presents numerical checks and examples in dimensions five and seven.

\subsection*{Acknowledgements}
The first author gratefully acknowledges financial support from Da Lat University through research project No.~B2026-TDL-01-KPT. The corresponding author gratefully acknowledges financial support from Ho Chi Minh City University of Agriculture and Forestry through the institutional level science and technology project, grant number CS-CB25-KH-02.

\section{The geometry of lines and conics}\label{sec:geometry}
By a conic, we mean a plane quadratic subscheme of $\PP^r$. Let $G_2=\Gr(3,r+1)$ be the Grassmannian of planes in $\PP^r$, with tautological subbundle $S$ of rank three. Put
\[
q:P=\PP(\Sym^2S^\vee)\longrightarrow G_2.
\]
We use the convention that $\PP(E)$ parametrizes one-dimensional subspaces of $E$. Thus, a point of $P$ is a pair $(\Pi,[\Phi])$, where $\Pi\subset\PP^r$ is a plane and $\Phi$ is a nonzero quadratic form on $\Pi$, up to scalar. Its zero scheme is a plane conic, possibly singular. Write
\[
\xi=c_1(\OO_P(1)).
\]
For $d\geq2$, multiplication by the universal quadratic form gives an exact sequence
\begin{equation}\label{eq:Em}
0\longrightarrow
q^*\Sym^{d-2}S^\vee\otimes\OO_P(-1)
\longrightarrow q^*\Sym^dS^\vee
\longrightarrow E_d
\longrightarrow0.
\end{equation}
At a point $(\Pi,[\Phi])$, the first map is multiplication by $\Phi$. It is injective because the polynomial ring of $\Pi$ is an integral domain and $\Phi\neq0$. This remains true when $\Phi$ defines a reducible conic or a double line. Hence, the cokernel $E_d$ is a vector bundle, with fiber
\[
(E_d)_{(\Pi,[\Phi])}
\simeq H^0(C,\OO_C(d)),
\quad C=V(\Phi)\subset\Pi.
\]
Its rank is
\[
\rank E_d
=
\binom{d+2}{2}-\binom{d}{2}
=2d+1.
\]
For $d=1$, set $E_1=q^*S^\vee$; its fiber is likewise $H^0(C,\OO_C(1))$, and its rank is three.

\begin{lemma}\label{lem:geometry}
Assume \eqref{eq:dimension}. For general defining equations of $X$, the induced section of
\[
E_X=\bigoplus_{i=1}^s E_{m_i}
\]
has a finite reduced zero scheme, disjoint from the discriminant divisor of singular plane conics in $P$. Thus, every conic on $X$ is smooth. The Fano scheme $F=F_1(X)$ is empty or smooth of pure dimension $h$.
\end{lemma}
\begin{proof}
For each $d\geq1$, restriction to planes gives a surjective bundle homomorphism
\[
H^0(\PP^r,\OO_{\PP^r}(d))\otimes\OO_{G_2}
\longrightarrow \Sym^dS^\vee.
\]
After pulling back to $P$ and composing with the quotient map defining $E_d$, we obtain a surjection
\[
H^0(\PP^r,\OO_{\PP^r}(d))\otimes\OO_P
\longrightarrow E_d.
\]
Consequently, the vector space
\[
W=\bigoplus_{i=1}^s
H^0(\PP^r,\OO_{\PP^r}(m_i))
\]
generates $E_X$.
Bertini's theorem for a generating space of sections shows that a general section has a smooth zero scheme of expected dimension zero, since
\[
\rank E_X
=\sum_{i=1}^s(2m_i+1)
=2m+s
=\dim P
=3r-1.
\]
As $P$ is projective, this zero scheme is finite and reduced. Let $\Delta\subset P$ be the discriminant divisor parametrizing singular plane conics. The incidence scheme
\[
I_\Delta
=
\{(z,w)\in\Delta\times W:w(z)=0\}
\]
is the total space of the kernel of the surjective evaluation map
\[
W\otimes\OO_\Delta\longrightarrow E_X|_\Delta.
\]
Hence
\[
\dim I_\Delta
=
\dim\Delta+\dim W-\rank E_X
=
\dim W-1.
\]
The projection $I_\Delta\to W$ is proper because $\Delta$ is projective. Its image is therefore closed and has dimension at most $\dim W-1$. A general tuple of defining equations consequently has no zero on $\Delta$. Similarly, $W$ generates the vector bundle
\[
A_X=\bigoplus_{i=1}^s\Sym^{m_i}U^\vee
\]
on $G_1$. The induced section has a zero scheme $F=F_1(X)$. Its expected dimension is
\[
\dim G_1-\rank A_X
=
2r-2-(m+s)
=
\frac{n-3}{2}
=h.
\]
Another application of Bertini's theorem shows that $F$ is empty or smooth of pure dimension $h$. The required properties hold simultaneously on a Zariski open dense subset of $W$.
\end{proof}

\begin{lemma}\label{lem:normal}
Under the hypotheses of Theorem~\ref{thm:main}, the normal bundle of every line $\ell\subset X$ splits as 
\begin{equation}\label{eq:splitting}
N_{\ell/X}
\simeq
\OO_\ell^{\oplus h}
\oplus
\OO_\ell(-1)^{\oplus(h+2)}.
\end{equation}
In particular,
\[
H^1(\ell,N_{\ell/X})=0,
\quad
H^0(\ell,N_{\ell/X}(-1))=0,
\]
and the homomorphism $\delta$ in \eqref{eq:delta} is injective on every fiber.
\end{lemma}
\begin{proof}
Since $X$ is a smooth complete intersection, the normal sequence for $\ell\subset X\subset\PP^r$ is
\begin{equation}\label{eq:normal}
0\longrightarrow N_{\ell/X}
\longrightarrow \OO_\ell(1)\otimes Q_\ell
\longrightarrow
\bigoplus_{i=1}^s\OO_\ell(m_i)
\longrightarrow0.
\end{equation}
The map on global sections on the right is the differential, at $[\ell]$, of the section of $A_X$ defining $F$. Since $F$ is smooth of the expected dimension, this differential is surjective. As
\[
H^1(\ell,\OO_\ell(1)\otimes Q_\ell)=0,
\]
the cohomology sequence gives
\[
H^1(\ell,N_{\ell/X})=0.
\]
We next prove that $H^0(\ell,N_{\ell/X}(-1))=0$. Suppose otherwise. Twisting \eqref{eq:normal} by $\OO_\ell(-1)$, a nonzero section determines a nonzero vector $v\in Q_\ell$ whose normal derivatives of all the defining equations vanish along $\ell$. Choose a lift of $v$ to $\CC^{r+1}$, and let $\Pi\simeq\PP^2$ be the plane spanned by this lift and $\ell$. Choose homogeneous coordinates $[x_0:x_1:z]$ on $\Pi$ such that $\ell=(z=0)$ and the $z$-direction corresponds to the chosen lift. For each defining equation $f_i$, the restriction $f_i|_\Pi$ has zero constant term in $z$, since $\ell\subset X$, and zero linear term, since the normal derivative in direction $v$ vanishes. Thus $f_i|_\Pi$ is divisible by $z^2$; if $m_i=1$, it is identically zero. Hence the planar double line
\[
V_\Pi(z^2)\subset\Pi
\]
is contained in $X$. This is a singular plane conic, contrary to Lemma~\ref{lem:geometry}. Therefore
\[
H^0(\ell,N_{\ell/X}(-1))=0.
\]
By the splitting theorem for vector bundles on $\PP^1$, write
\[
N_{\ell/X}
\simeq\bigoplus_{j=1}^{n-1}\OO_\ell(a_j).
\]
The vanishing of $H^1(\ell,N_{\ell/X})$ gives $a_j\geq-1$, whereas the vanishing of $H^0(\ell,N_{\ell/X}(-1))$ gives $a_j\leq0$. Thus every $a_j$ is either $0$ or $-1$. Moreover,
\[
\rank N_{\ell/X}=n-1=2h+2,
\quad
\deg N_{\ell/X}
=(r+1-m)-2=-h-2.
\]
It follows that precisely $h+2$ summands have degree $-1$ and the remaining $h$ have degree zero. This proves \eqref{eq:splitting}. Finally, twisting \eqref{eq:normal} by $\OO_\ell(-1)$
and taking global sections gives
\[
\ker\delta_\ell=H^0(\ell,N_{\ell/X}(-1))=0.
\]
\end{proof}

\begin{corollary}\label{cor:conicnormal}
Under the hypotheses of Theorem~\ref{thm:main}, every embedded conic $C\subset X$ has the normal bundle
\[
N_{C/X}\simeq
\OO_{\PP^1}(-1)^{\oplus(n-1)},
\]
where $C$ is identified with $\PP^1$. Here $\OO_{\PP^1}(1)$ has degree one, whereas the hyperplane bundle $\OO_C(1)$ corresponds to $\OO_{\PP^1}(2)$.
\end{corollary}

\begin{proof}
By Lemma~\ref{lem:geometry}, the Hilbert point of $C$ is isolated and reduced. Hence $H^0(C,N_{C/X})=0$. Write
\[
N_{C/X}\simeq
\bigoplus_{j=1}^{n-1}\OO_{\PP^1}(b_j).
\]
The vanishing implies $b_j\leq-1$ for every $j$. On the other hand, adjunction and \eqref{eq:dimension} give
\[
\sum_{j=1}^{n-1}b_j
=
\deg N_{C/X}
=
2(r+1-m)-2
=
-(n-1).
\]
Thus every $b_j=-1$.
\end{proof}

Let
\[
p:\LL=\PP(U|_F)\longrightarrow F
\]
be the universal line, and let $e:\LL\to X$ be the evaluation. The morphism $(e,p):\LL\hookrightarrow X\times F$ is a closed immersion. Write
\[
N=N_{\LL/(X\times F)}
\]
for its normal bundle, whose restriction to a fiber $\ell$ of $p$ is $N_{\ell/X}$. By \eqref{eq:splitting}, cohomology and base change give
\[
p_*N\simeq T_F,
\quad
R^1p_*N=0,
\]
where the first identification is the usual description of the tangent bundle of the smooth Fano scheme $F$. The evaluation homomorphism
\[
p^*p_*N=p^*T_F\longrightarrow N
\]
restricts on each fiber to the inclusion of the trivial summands in \eqref{eq:splitting}. It therefore has constant rank $h$, and its cokernel $\mathscr R$ is locally free, with
\[
\mathscr R|_\ell
\simeq\OO_\ell(-1)^{\oplus(h+2)}.
\]
Thus, $\mathscr R(1)$ is trivial on every fiber. Cohomology and base change show that
\[
B':=p_*(\mathscr R(1))
\]
is locally free of rank $h+2$, and evaluation gives
\[
\mathscr R\simeq\OO_{\LL}(-1)\otimes p^*B'.
\]
We identify $B'$ with the bundle $B$ defined in \eqref{eq:B}. The relative normal sequence is
\[
0\longrightarrow N
\longrightarrow \OO_{\LL}(1)\otimes p^*(Q|_F)
\longrightarrow
\bigoplus_{i=1}^s\OO_{\LL}(m_i)
\longrightarrow0.
\]
Twisting by $\OO_{\LL}(-1)$ and pushing forward gives
\[
0\longrightarrow Q|_F
\xrightarrow{\delta}
\bigoplus_{i=1}^s\Sym^{m_i-1}U^\vee|_F
\longrightarrow R^1p_*(N(-1))
\longrightarrow0.
\]
Here $p_*(N(-1))=0$ by \eqref{eq:splitting}. On the other hand, twisting the evaluation sequence gives
\[
0\longrightarrow
\OO_{\LL}(-1)\otimes p^*T_F
\longrightarrow N(-1)
\longrightarrow
\OO_{\LL}(-2)\otimes p^*B'
\longrightarrow0.
\]
Our convention for projectivization gives
\[
\omega_p
\simeq\OO_{\LL}(-2)\otimes p^*\det U^\vee,
\qquad
R^1p_*\OO_{\LL}(-2)\simeq\det U.
\]
Since both direct images of $\OO_{\LL}(-1)$ vanish, relative duality and the projection formula yield
\[
R^1p_*(N(-1))\simeq\det U\otimes B'.
\]
Consequently,
\[
B'
\simeq\det U^\vee\otimes\coker(\delta)
=B.
\]
We have, therefore, established the exact sequence
\begin{equation}\label{eq:evaluation}
0\longrightarrow p^*T_F
\longrightarrow N
\longrightarrow\OO_{\LL}(-1)\otimes p^*B
\longrightarrow0.
\end{equation}
Taking classes in $K^0(F)$ also gives \eqref{eq:Kclass}.

\section{The compactified double cover family}\label{sec:covers}
Let $\HH=\Mbar_{0,0}(\LL/F,2)$ be the relative stack of unmarked genus zero stable maps of degree two to the fibers of $p:\LL\to F$. An object over a scheme $S$ consists of a morphism $a_S:S\to F$, a family of connected nodal curves $\MC_S\to S$ of arithmetic genus zero, and a stable map
\[
f_S:\MC_S\longrightarrow\LL\times_F S
\]
of degree two on every geometric fiber. The morphism
\[
\rho:\HH\longrightarrow F
\]
records the target line. Its fiber over $[\ell]\in F$ is $\Mbar_{0,0}(\ell,2)\simeq\Mbar_{0,0}(\PP^1,2)$. At the nodal boundary, the source has two components, each mapped isomorphically onto the same reduced line $\ell$. Write
\[
\pi:\MC\longrightarrow\HH,
\quad
f:\MC\longrightarrow
\LL_{\HH}:=\LL\times_F\HH
\]
for the universal source curve and universal stable map. Let
\[
\tau:\LL_{\HH}\longrightarrow\LL,
\quad
p_{\HH}:\LL_{\HH}\longrightarrow\HH
\]
be the projections. By construction,
\[
p\circ\tau=\rho\circ p_{\HH},
\quad
p_{\HH}\circ f=\pi.
\]
We use the notation
\[
N_{\HH}=\tau^*N,
\quad
\OO_{\LL_{\HH}}(j)=\tau^*\OO_{\LL}(j).
\]
Thus, the universal stable map to $X$ is 
\[e\circ\tau\circ f:\MC\to X.\]

Put
\[
q_1:Z=\PP(\Sym^2U^\vee|_F)\longrightarrow F,
\quad H=c_1(\OO_Z(1)).
\]
A point of $Z$ consists of a line $\ell\subset X$ and a nonzero binary quadratic form on $\ell$, up to scalar. Its zero divisor has degree two, possibly supported at a single point. Let
\[
p_Z:\LL_Z:=\LL\times_F Z\longrightarrow Z
\]
be the projection. The tautological inclusion
\[
\OO_Z(-1)\longrightarrow
q_1^*(\Sym^2U^\vee|_F)
\]
determines a universal branch section
\[
b\in H^0\!\left(
\LL_Z,\OO_{\LL_Z}(2)\otimes p_Z^*\OO_Z(1)
\right).
\]
Let $\gamma:\mathcal G\to Z$ be the gerbe of square roots of $\OO_Z(-1)$, with universal line bundle $T$ and an isomorphism
\begin{equation}\label{eq:root}
T^{\otimes2}\simeq\gamma^*\OO_Z(-1).
\end{equation}
Write
\[
p_{\mathcal G}:\LL_{\mathcal G}:=
\LL\times_F\mathcal G\longrightarrow\mathcal G,
\quad
K=\OO_{\LL_{\mathcal G}}(-1)\otimes p_{\mathcal G}^*T.
\]
Pulling back $b$ gives a section of $K^{-2}$, or equivalently, a homomorphism
\[
K^{\otimes2}\longrightarrow\OO_{\LL_{\mathcal G}}.
\]
Together with the usual scalar multiplication, this defines a rank-two algebra
\begin{equation}\label{eq:algebra}
\mathscr A
=
\OO_{\LL_{\mathcal G}}\oplus K
=
\OO_{\LL_{\mathcal G}}\oplus
\bigl(\OO_{\LL_{\mathcal G}}(-1)
      \otimes p_{\mathcal G}^*T\bigr).
\end{equation}
The associated double cover is
\[
\Spec_{\LL_{\mathcal G}}\mathscr A
\longrightarrow\LL_{\mathcal G}.
\]
Locally, it has the equation $w^2=b$.

\begin{lemma}\label{lem:gerbe}
The branch construction induces an isomorphism
\[
\HH\simeq\mathcal G
\]
over $F$. Under this identification, the branch morphism $\gamma:\HH\to Z$ is a $\mu_2$-gerbe and $\rho=q_1\circ\gamma$. Consequently, $\HH$ is a smooth, proper Deligne-Mumford stack of pure dimension $h+2$, and $\rho$ is smooth of relative dimension two. Moreover,
\[
\gamma_*1=\frac12\,1_Z
\quad\text{in }A^0(Z)_{\mathbb Q}.
\]
\end{lemma}
\begin{proof}
We construct the equivalence on families over an arbitrary locally Noetherian $\CC$-scheme $S$. Let $a_S:S\to F$, and write
\[
p_S:L_S=\PP(a_S^*(U|_F))\longrightarrow S.
\]
An object of $\HH(S)$ is a family of genus zero stable maps
\[
f_S:C_S\longrightarrow L_S
\]
of degree two on every geometric fiber. We first show that no geometric fiber has a contracted component. The dual graph of a genus zero nodal curve is a tree. Since the total degree is two, at most two vertices have positive degree. Every contracted component must meet the other components in at least three points, by stability and the absence of markings. Its vertex therefore has valence at least three. A tree containing such a vertex has at least three leaves. Each leaf must have positive degree, since a contracted component meeting the rest of the curve at only one point would be unstable. This contradicts the total degree two. Consequently, each geometric fiber is either a smooth rational curve mapping with degree two to the target line, or two smooth rational components meeting at one node, each mapped isomorphically onto that line. The morphism $f_S$ is proper and has finite geometric fibers, hence is finite. On each geometric fiber over $S$, it is flat of degree two. For a smooth source, this follows because a finite morphism between smooth integral curves is flat. In the nodal case, near the image of the node, the cover has the local form
\[
k[z]\longrightarrow k[z,w]/(w^2-z^2),
\]
which is free of rank two over $k[z]$, with basis $1,w$. Away from the node, the two components map isomorphically onto the target. Since $C_S$ is flat over $S$, the fiberwise criterion for flatness shows that $f_S$ is flat. Thus
\[
\mathscr A_S:=f_{S*}\OO_{C_S}
\]
is a locally free $\OO_{L_S}$-algebra of rank two. The trace map satisfies $\tr(1)=2$. Since $2$ is invertible, the unit inclusion splits canonically by $\frac12\tr$, giving
\[
\mathscr A_S=\OO_{L_S}\oplus K,
\quad
K=\ker\bigl(\tr:\mathscr A_S\to\OO_{L_S}\bigr),
\]
where $K$ is invertible. Multiplication in $\mathscr A_S$ restricts to a homomorphism
\[
K^{\otimes2}\longrightarrow\OO_{L_S}.
\]
Indeed, locally choose a generator $w$ of $K$ and write $w^2=a+bw$. Relative to the basis $1,w$, multiplication by $w$ has trace $b$. Since $\tr(w)=0$, we have $b=0$; thus, $w^2=a$ lies in the scalar summand. For each geometric point $s\to S$, the finiteness of $f_s$ and the trace decomposition give
\[
1=\chi(\OO_{C_s})
=\chi(\OO_{L_s})+\chi(K_s)
=1+\chi(K_s).
\]
Thus $\chi(K_s)=0$. Since $L_s\simeq\PP^1$, the Riemann-Roch theorem gives $\deg K_s=-1$, so
\[
K_s\otimes\OO_{L_s}(1)\simeq\OO_{L_s}.
\]
In particular, this line bundle has one global section and vanishing first cohomology on every geometric fiber. Cohomology and base change therefore show that
\[
T_S:=p_{S*}(K\otimes\OO_{L_S}(1))
\]
is invertible and commutes with base change. The evaluation homomorphism
\[
p_S^*T_S\longrightarrow K\otimes\OO_{L_S}(1)
\]
is an isomorphism on every geometric fiber, hence an isomorphism. Consequently,
\[
K\simeq\OO_{L_S}(-1)\otimes p_S^*T_S.
\]
Under this identification, multiplication $K^{\otimes2}\to\OO_{L_S}$ becomes a homomorphism
\[
p_S^*(T_S^{\otimes2})\longrightarrow\OO_{L_S}(2).
\]
By adjunction and the projective bundle formula
\[
p_{S*}\OO_{L_S}(2)
\simeq a_S^*(\Sym^2U^\vee|_F),
\]
it is equivalent to
\[
T_S^{\otimes2}
\longrightarrow a_S^*(\Sym^2U^\vee|_F).
\]
This homomorphism is nonzero on every geometric fiber. Otherwise, the corresponding cover would have algebra $\OO_{L_s}\oplus K_s$ with square zero ideal $K_s$, contradicting the reducedness of the nodal curve $C_s$. It therefore identifies $T_S^{\otimes2}$ with a line subbundle of $a_S^*(\Sym^2U^\vee|_F)$. By the lines convention for projectivization, this subbundle determines a morphism
\[
g_S:S\longrightarrow Z,
\quad q_1\circ g_S=a_S,
\]
together with an isomorphism
\[
T_S^{\otimes2}\simeq g_S^*\OO_Z(-1).
\]
These are precisely the square root data defining an object of $\mathcal G(S)$. Conversely, let an object of $\mathcal G(S)$ be given by a morphism $g_S:S\to Z$, an invertible sheaf $T_S$, and an isomorphism
\[
T_S^{\otimes2}\simeq g_S^*\OO_Z(-1).
\]
Put $a_S=q_1\circ g_S$ and
\[
p_S:L_S=\PP(a_S^*(U|_F))\longrightarrow S.
\]
Pulling back the tautological inclusion on $Z$ gives
\[
T_S^{\otimes2}
\longrightarrow a_S^*(\Sym^2U^\vee|_F).
\]
By adjunction, this defines a homomorphism
\[
K_S^{\otimes2}\longrightarrow\OO_{L_S},
\quad
K_S=\OO_{L_S}(-1)\otimes p_S^*T_S.
\]
It determines a commutative algebra structure on
\[
\mathscr A_S=\OO_{L_S}\oplus K_S,
\]
with unit given by the first summand. Locally, this algebra has the form $\OO_{L_S}[w]/(w^2-b)$. Set
\[
C_S=\Spec_{L_S}\mathscr A_S.
\]
The resulting morphism $f_S:C_S\to L_S$ is finite flat of degree two. Since $L_S\to S$ is flat and proper, so is $C_S\to S$. On each geometric fiber, the branch section is a nonzero binary quadratic. If it has two distinct zeros, the cover is smooth and connected, and the Riemann-Hurwitz formula shows that its genus is zero. If it has a double zero, it is the square of a linear form, up to a nonzero scalar. The cover then consists of two copies of the target meeting transversely at that zero, with a local equation
\[
w^2=z^2.
\]
Thus, every geometric fiber is a connected nodal curve of arithmetic genus zero. In both cases, the map has degree two and no contracted components. Its automorphism group over the target has order two: it is generated by the deck involution in the smooth case and by the interchange of the two components in the nodal case. Hence, $f_S$ is a family of stable maps and defines an object of $\HH(S)$.

The trace decomposition commutes with arbitrary base change because the cover is finite locally free. The formation of $T_S$ and its evaluation isomorphism also commutes with base change, by the cohomology and base change argument above. The construction of the algebra from the square root data is likewise compatible with base change.

These constructions also identify morphisms. An isomorphism of covers induces an isomorphism of their algebras preserving the unit, trace, and multiplication. It therefore identifies the trace zero summands and induces an isomorphism of the corresponding line bundles $T_S$ compatible with their square root isomorphisms. Conversely, an isomorphism of square root data induces an algebra isomorphism and hence an isomorphism of covers. The two constructions are thus mutually inverse as functors of groupoids, compatibly with pullback.

This proves the equivalence over locally Noetherian bases, including nonreduced bases. Since the families, line bundles, and homomorphisms involved are of finite presentation, the equivalence extends to arbitrary $\CC$-schemes by Noetherian approximation and descent. Hence
\[
\HH\simeq\mathcal G
\]
as stacks over $F$. Under this identification, write
\[
\gamma:\HH\longrightarrow Z
\]
for the gerbe projection and $T$ for the universal square root line bundle. Then
\[
\rho=q_1\circ\gamma.
\]
Since $F$ is smooth and proper of pure dimension $h$, the projective bundle $Z$ is smooth and proper of pure dimension $h+2$. The morphism $\gamma$ is a $\mu_2$-gerbe, so $\HH$ is a smooth proper Deligne-Mumford stack of the same dimension, and $\rho$ is smooth of relative dimension two. Finally, integration along the gerbe gives
\[
\gamma_*1=\frac12\,1_Z
\quad\text{in }A^0(Z)_{\mathbb Q}.
\]
\end{proof}

The preceding construction gives an explicit description of the branch morphism on families. For the general construction of branch morphisms for stable maps, see~\cite{FantechiP}. Define
\begin{equation}\label{eq:A}
A:=R^1\pi_*f^*\OO_{\LL_{\HH}}(-1).
\end{equation}
We identify this sheaf using the algebra \eqref{eq:algebra}. Since $f$ is finite, the projection formula gives
\[
f_*f^*\OO_{\LL_{\HH}}(-1)
\simeq
\OO_{\LL_{\HH}}(-1)
\oplus
\bigl(\OO_{\LL_{\HH}}(-2)\otimes p_{\HH}^*T\bigr).
\]
Using $\pi=p_{\HH}\circ f$, base change, and
\[
R^1p_*\OO_{\LL}(-1)=0,
\quad
R^1p_*\OO_{\LL}(-2)\simeq\det U,
\]
we obtain
\[
A\simeq\rho^*\det U\otimes T.
\]
In particular, $A$ is a line bundle on the entire compactified space $\HH$. Set $\theta=c_1(A)$. Since
\[
c_1(\det U)=-\sigma,
\quad
2c_1(T)=-\gamma^*H
\]
by \eqref{eq:root}, we have
\begin{equation}\label{eq:theta}
\theta=-\rho^*\sigma-\frac12\gamma^*H
\quad\text{in }A^1(\HH)_{\mathbb Q}.
\end{equation}

\begin{proposition}\label{prop:obstruction}
Under the hypotheses of Theorem~\ref{thm:main}, composition with the inclusion of each line into $X$ induces an open and closed immersion
\[
i:\HH\longrightarrow\Mbar_{0,0}(X,2).
\]
Its image is smooth of pure dimension $h+2$. The complement is a finite reduced scheme whose points correspond to smooth embedded conics in $X$. The obstruction sheaf induced by the stable map perfect obstruction theory restricts to the vector bundle
\begin{equation}\label{eq:obstruction}
\Ob\simeq\rho^*B\otimes A,
\quad \rank\Ob=h+2,
\end{equation}
on the entire stack $\HH$, including its nodal boundary. The restriction of the virtual fundamental class is
\[
[\Mbar_{0,0}(X,2)]^{\mathrm{vir}}\big|_{\HH}
=
c_{h+2}(\Ob)\cap[\HH].
\]
Each point of the complement contributes one to $N_2(X)$.
\end{proposition}
\begin{proof}
Put
\[
\MM=\Mbar_{0,0}(X,2),
\quad
\YY=\Mbar_{0,0}(\PP^r,2).
\]
Let $[g:C\to X]$ be a geometric point of $\MM$. The dual graph argument in Lemma~\ref{lem:gerbe} shows that $C$ has no contracted components. If $C$ is irreducible, it is a smooth rational curve, and $g$ is either a double cover of a line or an isomorphism onto a smooth conic. If $C$ is reducible, it consists of two rational components meeting at one node, each mapped isomorphically onto a line. These image lines either coincide, giving a boundary point of $\HH$, or are distinct and intersect at the image of the node. The latter case would give a singular plane conic contained in $X$, contrary to Lemma~\ref{lem:geometry}. Thus, every geometric point of $\MM$ represents either a degree two map to a line, parametrized by $\HH$, or an embedded smooth conic. We now show that this separation holds scheme theoretically.

Let
\[
\pi_{\YY}:\MC_{\YY}\longrightarrow\YY,
\quad
g_{\YY}:\MC_{\YY}\longrightarrow\PP^r
\]
be the universal curve and universal stable map. On every geometric fiber $C$, the line bundle $g_{\YY}^*\OO_{\PP^r}(1)|_C$ has total degree two and nonnegative degree on each component. Since $C$ has arithmetic genus zero, its first cohomology vanishes, and the Riemann--Roch formula gives
\[
h^0(C,g_{\YY}^*\OO_{\PP^r}(1)|_C)=3.
\]
Cohomology and base change, therefore, imply
\[
R^1\pi_{\YY*}g_{\YY}^*\OO_{\PP^r}(1)=0,
\]
and
\[
\mathscr W
:=\pi_{\YY*}g_{\YY}^*\OO_{\PP^r}(1)
\]
is locally free of rank three and commutes with arbitrary base change. The pullback of linear forms defines a homomorphism
\[
\varepsilon:
H^0(\PP^r,\OO_{\PP^r}(1))\otimes\OO_{\YY}
\longrightarrow\mathscr W.
\]
At a geometric point, its rank is one more than the dimension of the linear span of the image curve. In particular, its rank is at least two since a degree two stable map has a nonconstant image. Let $D\subset\YY$ be the closed determinantal substack defined by the $3\times3$ minors of $\varepsilon$. On $D$, the homomorphism has rank exactly two at every geometric point. Locally, choose an invertible $2\times2$ minor. Elementary changes of bases put $\varepsilon|_D$ into the form
\[
\begin{pmatrix}
I_2 & 0\\
0   & R
\end{pmatrix}.
\]
Every entry of $R$ is then a $3\times3$ minor, up to sign, and hence vanishes on $D$. Thus, $\varepsilon|_D$ has image a locally free direct summand of $\mathscr W|_D$ of rank two. This argument holds over nonreduced bases as well. 

Let $S\to D$, and let
\[
g_S:C_S\longrightarrow\PP^r_S
\]
be the corresponding family of stable maps. Write $\mathscr E_S$ for the rank two image of the pullback homomorphism $\varepsilon$. The quotient
\[
H^0(\PP^r,\OO(1))\otimes\OO_S
\twoheadrightarrow\mathscr E_S
\]
determines a family of lines
\[
L_S=\PP(\mathscr E_S^\vee)\hookrightarrow\PP^r_S,
\]
where the dual is required by our line convention for projectivization. The evaluation homomorphism on $C_S$ factors through $\pi_S^*\mathscr E_S$. Since the original evaluation homomorphism is surjective, the induced homomorphism
\[
\pi_S^*\mathscr E_S
\longrightarrow g_S^*\OO_{\PP^r_S}(1)
\]
is also surjective. It therefore defines a factorization
\[
C_S\xrightarrow{f_S}L_S\hookrightarrow\PP^r_S
\]
of $g_S$. The family $L_S$ is uniquely determined by the image of $\varepsilon_S$.

Conversely, suppose that a family of degree-two stable maps factors through a family of lines $p_S:L_S\to S$. Pullback of sections gives
\[
p_{S*}\OO_{L_S}(1)
\longrightarrow
\pi_{S*}g_S^*\OO_{\PP^r_S}(1).
\]
The source is locally free of rank two, and this map is injective on every geometric fiber because each cover surjects onto its target line. Its image is therefore a rank two subbundle. Since linear forms on $\PP^r_S$ generate $p_{S*}\OO_{L_S}(1)$, the evaluation homomorphism factors through this subbundle. Hence, all its $3\times3$ minors vanish, and the corresponding morphism $S\to\YY$ factors through $D$. These constructions identify $D$, including its scheme structure, with the stack of degree-two maps to families of lines. Now let $S\to D\times_{\YY}\MM$. By the proof of Lemma~\ref{lem:gerbe}, the induced map
\[
f_S:C_S\longrightarrow L_S
\]
is finite faithfully flat. Its composite with $L_S\hookrightarrow\PP^r_S$ factors through $X_S:=X\times S$. Thus, the ideal of $X_S$ restricts to zero after
pullback to $C_S$. Faithful flatness implies that
it already restricts to zero on $L_S$, so
\[
L_S\subset X_S.
\]
Therefore, the family of lines defines a morphism $S\to F$, and the cover defines an object of $\HH(S)$. Conversely, every object of $\HH(S)$ gives such a factorization. These identifications are compatible with morphisms and base change, yielding
\[
D\times_{\YY}\MM\simeq\HH.
\]
Since $D\hookrightarrow\YY$ is a closed immersion,
its base change
\[
i:\HH\longrightarrow\MM
\]
is a closed immersion as well.

We show that the closed immersion $i$ is also open. Since the stacks are of finite type over $\CC$, it suffices to check this at closed geometric points. Fix such a point $[f_0:C_0\to\ell]$ of $\HH$. Write $j_\ell:\ell\hookrightarrow X$ for the inclusion, and identify its image in $\MM$ with $[j_\ell\circ f_0]$.

The trace decomposition gives
\[
(f_0)_*\OO_{C_0}
\simeq\OO_\ell\oplus K_0,
\quad K_0\simeq\OO_\ell(-1).
\]
Since $f_0$ is finite, the projection formula yields
\[
H^1(C_0,f_0^*T_\ell)
\simeq
H^1(\ell,T_\ell)
\oplus H^1(\ell,T_\ell\otimes K_0)
=0.
\]
Here $T_\ell\simeq\OO_\ell(2)$, so the two summands are the first cohomology groups of $\OO_\ell(2)$ and $\OO_\ell(1)$. Similarly,
\[
H^0(C_0,f_0^*N_{\ell/X})
\simeq
H^0(\ell,N_{\ell/X})
\oplus H^0(\ell,N_{\ell/X}\otimes K_0).
\]
The second summand vanishes by \eqref{eq:splitting}. Thus pullback induces a canonical isomorphism
\[
H^0(\ell,N_{\ell/X})
\xrightarrow{\sim}
H^0(C_0,f_0^*N_{\ell/X}),
\]
and the source is naturally $T_{F,[\ell]}$. These arguments apply to both smooth and nodal sources. The deformation sequence associated with $\ell\hookrightarrow X$ gives an exact sequence
\[
0\longrightarrow
T_{[f_0]}\Mbar_{0,0}(\ell,2)
\longrightarrow
T_{[j_\ell\circ f_0]}\MM
\longrightarrow
H^0(C_0,f_0^*N_{\ell/X}),
\]
where the last map records the normal deformation of the stable map. We verify surjectivity of the last map. By Lemma~\ref{lem:gerbe}, $\rho:\HH\to F$ is smooth, so every tangent vector in $T_{F,[\ell]}$ lifts to $T_{[f_0]}\HH$. Under composition with $i$, its normal deformation is the pullback of the corresponding section of $N_{\ell/X}$. The preceding isomorphism therefore proves surjectivity. The tangent sequence of $\rho$ is
\[
0\longrightarrow
T_{[f_0]}\Mbar_{0,0}(\ell,2)
\longrightarrow
T_{[f_0]}\HH
\longrightarrow
T_{F,[\ell]}
\longrightarrow0.
\]
The differential of $i$ identifies the left terms of these sequences, while the right terms are identified by pullback. Hence
\[
di:T_{[f_0]}\HH
\xrightarrow{\sim}
T_{[j_\ell\circ f_0]}\MM.
\]

Choose an étale chart $V\to\MM$ and a closed point $v\in V$ above $[j_\ell\circ f_0]$.
The pullback $V_{\HH}=V\times_{\MM}\HH$ is a closed subspace of $V$. Let
\[
R=\OO_{V,v},
\qquad
R/I=\OO_{V_{\HH},v}.
\]
Since $V_{\HH}$ is étale over the smooth stack $\HH$, the ring $R/I$ is regular of dimension $d=h+2$. Equality of tangent spaces gives $\edim R=d$, so
\[
d=\dim(R/I)\leq\dim R\leq\edim R=d.
\]
Consequently, $R$ is regular. A regular local ring is an integral domain, and a nonzero ideal has quotient of strictly smaller dimension. Thus $I=0$. The ideal sheaf defining $V_{\HH}$ is coherent, so it vanishes on a neighborhood of $v$. Therefore, the closed immersion is an isomorphism onto a neighborhood at every closed point of $\HH$. It follows that $i$ is open as well as closed. In particular, $\MM$ is smooth along $\HH$.

The relative tangent sequence on the universal line is
\[
0\longrightarrow T_{\LL/F}
\longrightarrow e^*T_X
\longrightarrow N
\longrightarrow0.
\]
Write
\[
g=e\circ\tau\circ f:\MC\longrightarrow X
\]
for the universal stable map to $X$ over $\HH$. Pulling back the sequence gives
\[
0\longrightarrow f^*\tau^*T_{\LL/F}
\longrightarrow g^*T_X
\longrightarrow f^*N_{\HH}
\longrightarrow0.
\]
On every geometric fiber, the first term restricts to $f_0^*T_\ell$, whose first cohomology vanishes by the preceding argument. Cohomology and base change therefore give
\[
R^1\pi_*f^*\tau^*T_{\LL/F}=0.
\]
Since $\pi$ has one-dimensional fibers, the long exact sequence yields
\[
R^1\pi_*g^*T_X
\simeq R^1\pi_*f^*N_{\HH}.
\]

To identify the absolute obstruction sheaf, we must also account for the deformations of the source curve. At a geometric point $[f_0:C_0\to\ell]$, put
\[
g_0=j_\ell\circ f_0:C_0\longrightarrow X.
\]
The stable map deformation sequence contains
\[
\Ext^1(\Omega_{C_0},\OO_{C_0})
\longrightarrow H^1(C_0,g_0^*T_X)
\longrightarrow \Ob_{[g_0]}
\longrightarrow0.
\]
By naturality for $j_\ell:\ell\hookrightarrow X$, the first homomorphism factors as
\[
\Ext^1(\Omega_{C_0},\OO_{C_0})
\longrightarrow H^1(C_0,f_0^*T_\ell)
\longrightarrow H^1(C_0,g_0^*T_X).
\]
Its middle term is zero, so source deformations do not remove any obstruction classes.

This argument includes nodal sources. In particular, the local family
\[
w^2=z^2+t
\]
smooths the source node while keeping the target coordinate $z$ fixed. More generally, the factorization
\[
\MC\xrightarrow{f}\LL_{\HH}\longrightarrow X\times\HH
\]
induces the same factorization of the universal deformation morphisms. The homomorphism from source deformations to $R^1\pi_*g^*T_X$ factors through $R^1\pi_*f^*\tau^*T_{\LL/F}=0$. Thus, the identification holds at the level of sheaves, not merely at geometric points, and
\[
\Ob\simeq R^1\pi_*g^*T_X
\simeq R^1\pi_*f^*N_{\HH}
\]
on the entire compactified stack $\HH$.

Pulling back \eqref{eq:evaluation} by $\tau$ and $f$, and using $p\circ\tau\circ f=\rho\circ\pi$, gives
\[
0\longrightarrow\pi^*\rho^*T_F
\longrightarrow f^*N_{\HH}
\longrightarrow
f^*\OO_{\LL_{\HH}}(-1)\otimes\pi^*\rho^*B
\longrightarrow0.
\]
The fibers of $\pi$ are connected nodal curves of arithmetic genus zero. Hence
\[
\pi_*\OO_{\MC}\simeq\OO_{\HH},
\quad
R^1\pi_*\OO_{\MC}=0.
\]
By the projection formula,
\[
R^1\pi_*(\pi^*\rho^*T_F)=0.
\]
Moreover, $R^2\pi_*(\pi^*\rho^*T_F)=0$, since $\pi$ has one-dimensional fibers. The long exact sequence therefore yields
\[
\begin{aligned}
\Ob
&\simeq R^1\pi_*f^*N_{\HH}\\
&\simeq
R^1\pi_*\bigl(
f^*\OO_{\LL_{\HH}}(-1)\otimes\pi^*\rho^*B
\bigr)\\
&\simeq
\rho^*B\otimes R^1\pi_*f^*\OO_{\LL_{\HH}}(-1)\\
&=\rho^*B\otimes A.
\end{aligned}
\]
By \eqref{eq:A}, $A$ is a line bundle on the entire compactified stack $\HH$. Its local freeness can also be checked directly: on a smooth source, $f_0^*\OO_\ell(-1)\simeq\OO_{\PP^1}(-2)$; on a nodal source, it has degree $-1$ on each component. In either case,
\[
h^0(C_0,f_0^*\OO_\ell(-1))=0,
\quad
h^1(C_0,f_0^*\OO_\ell(-1))=1,
\]
so cohomology and base change apply. Consequently, $\Ob$ is locally free of rank $h+2$ across the nodal boundary as well. This proves \eqref{eq:obstruction}.

The stack $\YY$ is smooth because $\PP^r$ is convex. For each integer $d\geq1$, define
\[
V_2^{(d)}
=
\pi_{\YY*}g_{\YY}^*\OO_{\PP^r}(d).
\]
The subscript $2$ records the degree of the stable maps, while the superscript $d$ records the twist $\OO_{\PP^r}(d)$. The pullback line bundle has nonnegative degree on every component of every source curve. Hence
\[
R^1\pi_{\YY*}g_{\YY}^*\OO_{\PP^r}(d)=0,
\]
and $V_2^{(d)}$ is locally free of rank $2d+1$.
Put
\[
\mathcal V_X=\bigoplus_{i=1}^sV_2^{(m_i)}.
\]
The defining equations of $X$ induce a section $s_X$ of $\mathcal V_X$. Its zero locus is precisely $\MM$, including its scheme structure, since a family of maps factors through $X$ exactly when all the pullback defining equations vanish.

The resulting zero locus perfect obstruction theory agrees with the standard stable map obstruction theory for $X$. Since $\HH$ is a smooth open and closed substack of $\MM$, the differential of $s_X$ has constant rank along $\HH$ and gives an exact sequence
\[
0\longrightarrow T_{\HH}
\longrightarrow T_{\YY}|_{\HH}
\xrightarrow{\,ds_X\,}\mathcal V_X|_{\HH}
\longrightarrow\Ob
\longrightarrow0.
\]
The excess intersection formula~\cite{FultonIT}, or equivalently the induced perfect obstruction theory~\cite{BF}, therefore gives
\[
[\MM]^{\mathrm{vir}}|_{\HH}
=
c_{h+2}(\Ob)\cap[\HH]
=
c_{h+2}(\rho^*B\otimes A)\cap[\HH].
\]
It remains to consider the embedded smooth conics. For such a conic $C\subset X$, Lemma~\ref{lem:geometry} implies
\[
H^0(C,N_{C/X})=0.
\]
Since $C\simeq\PP^1$, its normal bundle has rank $n-1$ and degree
\[
\deg N_{C/X}
=
\deg(T_X|_C)-\deg T_C
=
2(r+1-m)-2.
\]
The Riemann-Roch formula and \eqref{eq:dimension} yield
\[
\chi(N_{C/X})
=
(n-1)+2(r+1-m)-2
=
n-3+2(r+1-m)
=0.
\]
Thus $H^1(C,N_{C/X})=0$ as well. The corresponding stable map point is reduced and unobstructed. Its automorphism group is trivial because the map is an embedding, so its contribution to $N_2(X)$ is one.
\end{proof}

For $j\geq0$, let $h_j(Y_1)$ denote the complete symmetric polynomial of degree $j$ in $Y_1=\{2x,x+y,2y\}$, defined by
\[
\sum_{j\geq0}h_j(Y_1)z^j
=\prod_{\eta\in Y_1}(1-\eta z)^{-1}.
\]
In particular, $h_0(Y_1)=1$. The subscripted notation $h_j$ is distinct from the dimension parameter $h=(n-3)/2$.

\begin{lemma}\label{lem:push}
For every integer $k\geq0$,
\begin{equation}\label{eq:push}
\rho_*(\theta^{2+k})
=
\frac{(-1)^k}{2^{k+3}}h_k(Y_1)
\quad\text{in }A^k(F)_{\mathbb Q}.
\end{equation}
Moreover,
\[
\rho_*1=\rho_*\theta=0.
\]
\end{lemma}
\begin{proof}
Recall that $\rho=q_1\circ\gamma$. With our convention for projectivization, the relation on $Z$ is
\[
(H+2x)(H+x+y)(H+2y)=0,
\]
where pullbacks from $F$ are understood. Set
\[
H'=H+2q_1^*\sigma,
\quad \sigma=x+y.
\]
Then
\[
(H'-2x)(H'-(x+y))(H'-2y)=0.
\]
Equivalently,
\[
H'^3-e_1(Y_1)H'^2+e_2(Y_1)H'-e_3(Y_1)=0,
\]
where $e_j(Y_1)$ denotes the elementary symmetric polynomial of degree $j$ in $Y_1$.

Put
\[
a_k=(q_1)_*(H'^{2+k})\quad(k\geq0),
\quad a_k=0\quad(k<0).
\]
The projective bundle pushforward formula gives
\[
(q_1)_*1=(q_1)_*H'=0,
\quad a_0=1.
\]
Multiplying the cubic relation by $H'^{k-1}$ and applying the projection formula yields
\[
a_k=e_1(Y_1)a_{k-1}
-e_2(Y_1)a_{k-2}
+e_3(Y_1)a_{k-3}
\quad(k\geq1).
\]
Thus
\[
\sum_{k\geq0}a_kz^k
=
\frac{1}{1-e_1(Y_1)z+e_2(Y_1)z^2-e_3(Y_1)z^3}
=
\prod_{\eta\in Y_1}(1-\eta z)^{-1}.
\]
By the definition of $h_k(Y_1)$, we obtain
\[
(q_1)_*(H'^{2+k})=h_k(Y_1).
\]

By \eqref{eq:theta}, $\theta=-\gamma^*H'/2$. Using $\gamma_*1=\frac12\,1_Z$ and the projection formula, we conclude that
\[
\begin{aligned}
\rho_*(\theta^{2+k})
&=
\left(-\frac12\right)^{2+k}
(q_1)_*\gamma_*\gamma^*(H'^{2+k})\\
&=
\frac12\left(-\frac12\right)^{2+k}h_k(Y_1)\\
&=
\frac{(-1)^k}{2^{k+3}}h_k(Y_1).
\end{aligned}
\]
Finally, $\rho_*1=\rho_*\theta=0$ because $\rho$ has relative dimension two.
\end{proof}

\begin{proof}[Proof of Theorem~\ref{thm:main}]
The geometric assertions follow from Lemmas~\ref{lem:geometry} and~\ref{lem:normal}. By Proposition~\ref{prop:obstruction}, every embedded smooth conic contributes one, while the double cover component contributes the top Chern class of its obstruction bundle. Hence
\[
N_2(X)=n_2(X)+D_2(X),
\quad
D_2(X)=\int_{\HH}c_{h+2}(\rho^*B\otimes A).
\]
Since $\rank B=h+2$ and $\theta=c_1(A)$,
\[
c_{h+2}(\rho^*B\otimes A)
=
\sum_{j=0}^{h+2}\rho^*c_j(B)\theta^{h+2-j}.
\]
Applying the projection formula and Lemma~\ref{lem:push}, we obtain
\[
\rho_*c_{h+2}(\rho^*B\otimes A)
=
\frac18\sum_{j=0}^h
c_j(B)\frac{(-1)^{h-j}}{2^{h-j}}h_{h-j}(Y_1).
\]
The terms with $j=h+1,h+2$ vanish because $\rho_*\theta=\rho_*1=0$. The generating function identity
\[
\sum_{k\geq0}(-1)^kh_k(Y_1)(t/2)^k
=\prod_{\eta\in Y_1}(1+\eta t/2)^{-1}
=c_{t/2}(\Sym^2U^\vee)^{-1}
\]
therefore gives
\[
D_2(X)
=
\frac18\int_F[t^h]\,
\frac{c_t(B)}{c_{t/2}(\Sym^2U^\vee)}.
\]
This proves \eqref{eq:main}. When $h=0$, the coefficient of $t^0$ is one, so $D_2(X)=\#F/8$, proving \eqref{eq:isolated}.
\end{proof}

\begin{proof}[Proof of Theorem~\ref{thm:coefficient}]
Let $\iota:F\hookrightarrow G_1$ be the inclusion. Since $F$ is the regular zero locus of the section of $A_X$ defined by the equations of $X$,
\[
\iota_*[F]=c_{m+s}(A_X)\cap[G_1].
\]
By \eqref{eq:Kclass},
\[
c_t(B)
=
\iota^*\left(
\frac{\prod_{i=1}^s
c_t(\det U^\vee\otimes\Sym^{m_i-1}U^\vee)}
{c_t(\det U^\vee\otimes Q)}
\right).
\]
Substituting this into \eqref{eq:main} and applying the projection formula gives \eqref{eq:ambient}. We next identify the coefficient class explicitly. Tensoring the tautological sequence by $\det U^\vee$ gives
\[
0\longrightarrow\det U^\vee\otimes U
\longrightarrow(\det U^\vee)^{\oplus(r+1)}
\longrightarrow\det U^\vee\otimes Q
\longrightarrow0.
\]
Since $U$ has rank two, $\det U^\vee\otimes U\simeq U^\vee$. Writing $\sigma=x+y$, we obtain
\[
c_t(\det U^\vee\otimes Q)
=
\frac{(1+\sigma t)^{r+1}}{(1+xt)(1+yt)}.
\]
Moreover,
\[
c_{t/2}(\Sym^2U^\vee)
=
(1+xt)(1+\sigma t/2)(1+yt).
\]
Consequently,
\[
c_t(\det U^\vee\otimes Q)\,
c_{t/2}(\Sym^2U^\vee)
=
(1+\sigma t)^{r+1}(1+\sigma t/2).
\]
The Chern roots of $\det U^\vee\otimes\Sym^{m_i-1}U^\vee$ are
\[
(j+1)x+(m_i-j)y,
\quad 0\leq j\leq m_i-1.
\]
Thus the coefficient class in \eqref{eq:ambient} is represented by the polynomial $K_X(x,y)$ in \eqref{eq:Kpoly}. Finally,
\[
c_{m+s}(A_X)
=
\prod_{i=1}^s\prod_{j=0}^{m_i}
\bigl(jx+(m_i-j)y\bigr).
\]
The product $c_{m+s}(A_X)K_X(x,y)$ has degree
\[
(m+s)+h=2r-2=\dim G_1.
\]
Applying the Grassmannian coefficient formula
\[
\int_{G_1}\varphi(x,y)
=
[x^ry^{r-1}]\bigl((x-y)\varphi(x,y)\bigr)
\]
to this symmetric polynomial proves \eqref{eq:ranktwo}.
\end{proof}

\begin{remark}\label{rem:JNS}
Put
\[
\mathcal E_X
=
\bigoplus_{i=1}^s\Sym^{m_i-1}U^\vee.
\]
Since $\rank\mathcal E_X=m$ and $m-r-1=h$, formula \eqref{eq:Kpoly} gives
\[
K_X(x,y)
=
[t^h]\,
\frac{(1+\sigma t)^h
c_{t/(1+\sigma t)}(\mathcal E_X)}
{1+\sigma t/2}.
\]
For every $k\geq0$, the binomial theorem yields
\[
[t^k]\frac{(1+\sigma t)^k}{1+\sigma t/2}
=
\left(\frac{\sigma}{2}\right)^k.
\]
Expanding in the Chern classes of $\mathcal E_X$, we therefore obtain
\[
K_X(x,y)
=
[t^h]\frac{c_t(\mathcal E_X)}{1-\sigma t/2}.
\]
Consequently,
\[
D_2(X)
=
\frac18\int_F[t^h]\,
\frac{\prod_{i=1}^s c_t(\Sym^{m_i-1}U^\vee)}
{1-\sigma t/2}.
\]
For a hypersurface, this recovers the correction formula of Jinzenji--Nakamura--Suzuki \cite[Theorem 5.1]{JNS}.
\end{remark}

\begin{corollary}\label{cor:denominators}
Under the hypotheses of Theorem~\ref{thm:main}, we have
\[
2^{h+3}D_2(X)\in\ZZ,\quad 2^{h+3}N_2(X)\in\ZZ.
\]
\end{corollary}
\begin{proof}
Expanding the coefficient in \eqref{eq:main} gives
\[
[t^h]\frac{c_t(B)}{c_{t/2}(\Sym^2U^\vee)}
=
\sum_{j=0}^h
\frac{(-1)^{h-j}}{2^{h-j}}
c_j(B)h_{h-j}(Y_1).
\]
Each product $c_j(B)h_{h-j}(Y_1)$ is an integral Chow class of codimension $h$ on $F$.
Consequently,
\[
2^{h+3}D_2(X)
=
\int_F\sum_{j=0}^h
(-1)^{h-j}2^j c_j(B)h_{h-j}(Y_1)
\in\ZZ.
\]
Here, the integral is an integer because $F$ is projective of pure dimension $h$. Finally, $n_2(X)\in\ZZ$ and $N_2(X)=n_2(X)+D_2(X)$, so
\[
2^{h+3}N_2(X)\in\ZZ.
\]
\end{proof}

\begin{corollary}\label{cor:five}
Under the hypotheses of Theorem~\ref{thm:main}, assume $n=5$. Then $F$ is empty or a smooth projective curve, and
\begin{equation}\label{eq:five}
D_2(X)
=
\frac18\int_F
\left(c_1(B)-\frac32\sigma\right).
\end{equation}
Moreover,
\[
c_1(B)
=
\left(
2+\frac12\sum_{i=1}^s m_i(m_i-1)
\right)\sigma.
\]
Consequently,
\begin{equation}\label{eq:fivefactor}
D_2(X)
=
\frac{\displaystyle 1+\sum_{i=1}^s m_i(m_i-1)}{16}
\int_F\sigma,
\end{equation}
where $\int_F\sigma$ is the degree of $F$ in its Pl\"ucker embedding and is understood to be zero when $F$ is empty.
\end{corollary}
\begin{proof}
Since $n=5$, we have $h=1$. Taking the coefficient of $t$ in \eqref{eq:main} gives
\[
D_2(X)
=
\frac18\int_F
\left(c_1(B)-\frac12c_1(\Sym^2U^\vee)\right).
\]
The Chern roots of $\Sym^2U^\vee$ are $2x,x+y,2y$, so $c_1(\Sym^2U^\vee)=3\sigma$. This proves \eqref{eq:five}. To compute $c_1(B)$, note that
\[
\rank\Sym^{m_i-1}U^\vee=m_i,
\quad
c_1(\Sym^{m_i-1}U^\vee)
=\frac{m_i(m_i-1)}2\,\sigma.
\]
Hence
\[
c_1(\det U^\vee\otimes\Sym^{m_i-1}U^\vee)
=
\left(m_i+\frac{m_i(m_i-1)}2\right)\sigma.
\]
The tautological sequence gives $\rank Q=r-1$ and $c_1(Q)=\sigma$, whence
\[
c_1(\det U^\vee\otimes Q)=r\sigma.
\]
Using \eqref{eq:Kclass} and $\displaystyle m=\sum_{i=1}^s m_i$, we obtain
\[
c_1(B)
=
\left(
m-r+\frac12\sum_{i=1}^s m_i(m_i-1)
\right)\sigma.
\]
Finally, \eqref{eq:dimension} and $n=r-s=5$ imply $m-r=2$. Thus
\[
c_1(B)
=
\left(
2+\frac12\sum_{i=1}^s m_i(m_i-1)
\right)\sigma.
\]
Substitution into \eqref{eq:five} yields
\[
D_2(X)
=
\frac{\displaystyle 1+\sum_{i=1}^s m_i(m_i-1)}{16}
\int_F\sigma,
\]
as claimed.
\end{proof}

\section{The localization formula}\label{sec:localization}
This section gives a localization formula for $N_2(X)$. Retain the bundles $V_2^{(d)}$ and $\mathcal V_X$ defined in the proof of Proposition~\ref{prop:obstruction}, and put
\[
E=\bigoplus_{i=1}^s\OO_{\PP^r}(m_i).
\]
For every $d\geq1$, the line bundle $g_{\YY}^*\OO_{\PP^r}(d)$ has total degree $2d$ and nonnegative degree on each component of every geometric source curve. Hence
\[
R^1\pi_{\YY*}g_{\YY}^*\OO_{\PP^r}(d)=0,
\quad
V_2^{(d)}=\pi_{\YY*}g_{\YY}^*\OO_{\PP^r}(d)
\]
is locally free of rank $2d+1$. Consequently,
\[
\mathcal V_X
=\pi_{\YY*}g_{\YY}^*E
=\bigoplus_{i=1}^sV_2^{(m_i)}
\]
is locally free of rank $2m+s$.

The defining equations of $X$ give the section $s_X$ of $\mathcal V_X$ with zero locus $\MM$ considered in Proposition~\ref{prop:obstruction}. Writing $j:\MM\hookrightarrow\YY$ for the inclusion, convex quantum Lefschetz gives
\[
j_*[\MM]^{\mathrm{vir}}
=c_{2m+s}(\mathcal V_X)\cap[\YY].
\]
Consequently,
\begin{equation}\label{eq:lefschetz}
N_2(X)
=
\int_{\YY}
\prod_{i=1}^s c_{2m_i+1}(V_2^{(m_i)}).
\end{equation}
Here $\YY=\Mbar_{0,0}(\PP^r,2)$ is smooth of dimension $3r-1=2m+s$ by \eqref{eq:dimension}, so the integrand has the required codimension.

Let $\mathbb T$ be a diagonal torus acting on $\PP^r$, and let $\lambda_i$ be the weight of $\OO_{\PP^r}(1)$ at the fixed point $p_i$, for $0\leq i\leq r$. Put
\[
P_i=\prod_{\substack{0\leq b\leq r\\b\ne i}}
(\lambda_i-\lambda_b).
\]
For an edge $e=(v,w)$ of degree $d_e$, with endpoint labels $i_v,i_w$, define
\[
\omega(v,e)=\frac{\lambda_{i_v}-\lambda_{i_w}}{d_e}.
\]
For a flag $F'=(v,e)$, write $\omega_{F'}=\omega(v,e)$. For distinct labels $i,j$ and an integer $d\geq1$, define
\begin{equation}\label{eq:edge}
E_{ij}(d)
=
\frac{(-1)^d(d!)^2}{d^{2d}}(\lambda_i-\lambda_j)^{2d}
\prod_{\substack{0\leq b\leq r\\b\ne i,j}}
\prod_{u=0}^d
\left(
\frac{u\lambda_i+(d-u)\lambda_j}{d}-\lambda_b
\right).
\end{equation}

\begin{proposition}\label{prop:local}
Let $\Gamma$ be a torus fixed graph of an unmarked genus zero stable map of degree two to $\PP^r$. Every vertex has valence one or two, and the corresponding fixed stratum is zero-dimensional. For each vertex $v$, put
\begin{equation}\label{eq:vertex}
V_v
=
P_{i_v}^{1-\operatorname{val}(v)}
\left(
\sum_{F'\ni v}\omega_{F'}^{-1}
\right)^{3-\operatorname{val}(v)}
\prod_{F'\ni v}\omega_{F'}.
\end{equation}
The sums and products run over the flags incident to $v$. If $N_\Gamma$ denotes the normal bundle of the fixed stratum, then
\begin{equation}\label{eq:euler}
e_{\mathbb T}(N_\Gamma)
=
\prod_v V_v
\prod_{e=(v,w)}E_{i_vi_w}(d_e).
\end{equation}
The localization contribution is divided by
\[
a_\Gamma
=
|\Aut(\Gamma)|\prod_e d_e,
\]
where $\Aut(\Gamma)$ consists of graph automorphisms preserving the vertex labels and edge degrees.
\end{proposition}
\begin{proof}
The edge and vertex factors are those in Kontsevich's localization formula~\cite{Kontsevich}; see also~\cite[Section 3]{Hiep2016}. Since the total degree is two, the graph has either one edge of degree two or two edges of degree one. Thus every vertex has valence one or two, and there are no contracted components or vertex moduli. For a valence one vertex with flag weight $\alpha$, \eqref{eq:vertex} gives $V_v=\alpha^{-1}$. For a valence two vertex with flag weights $\alpha,\beta$, it gives
\[
V_v
=
P_{i_v}^{-1}(\alpha^{-1}+\beta^{-1})\alpha\beta
=
\frac{\alpha+\beta}{P_{i_v}}.
\]
Together with the edge factors \eqref{eq:edge}, these yield \eqref{eq:euler}. The edge deck groups have orders $d_e$, and the label preserving graph automorphisms contribute $|\Aut(\Gamma)|$, giving the factor $a_\Gamma$.
\end{proof}

\begin{remark}
Proposition~\ref{prop:local} applies only in degree two. For general graphs, stable vertices contribute factors involving the cotangent line classes on $\Mbar_{0,\operatorname{val}(v)}$. These classes must be retained until integration over the corresponding vertex moduli spaces.
\end{remark}

For integers $k,d\geq1$, put
\[
S_k(i,j;d)=
\prod_{u=0}^{kd}
\frac{u\lambda_i+(kd-u)\lambda_j}{d}.
\]
The normalization sequence at the source nodes gives
\begin{equation}\label{eq:twist}
e_{\mathbb T}(V_2^{(k)}|_\Gamma)
=
\prod_{e=(v,w)}S_k(i_v,i_w;d_e)
\prod_v(k\lambda_{i_v})^{1-\operatorname{val}(v)}.
\end{equation}
At a valence two vertex, the last factor removes the duplicated contribution of the fiber at the node.

For a degree two edge joining vertices $v_1,v_2$ with labels $i<j$, put
\[
\begin{aligned}
S_X^{(2)}(i,j)
&=\prod_{a=1}^s S_{m_a}(i,j;2),\\
N^{(2)}(i,j)
&=V_{v_1}V_{v_2}E_{ij}(2),\\
C_{ij}^{(2)}
&=\frac{S_X^{(2)}(i,j)}{2N^{(2)}(i,j)}.
\end{aligned}
\]
The factor two is the order of the deck group.

For a chain with vertices $v_1,v_2,v_3$ and labels $(i,j,k)$ satisfying $i\ne j$ and $k\ne j$, put
\[
\begin{aligned}
S_X^{(1,1)}(i,j,k)
&=\prod_{a=1}^s
\frac{S_{m_a}(i,j;1)S_{m_a}(j,k;1)}{m_a\lambda_j},\\
N^{(1,1)}(i,j,k)
&=V_{v_1}V_{v_2}V_{v_3}E_{ij}(1)E_{jk}(1),\\
C_{ijk}^{(1,1)}
&=\frac{S_X^{(1,1)}(i,j,k)}{N^{(1,1)}(i,j,k)}.
\end{aligned}
\]
Each denominator $m_a\lambda_j$ cancels a common endpoint factor in the numerator. The vertex factors are evaluated using the incident flags of the indicated graph; in particular, $v_1$ and $v_3$ remain distinct when $i=k$. Reversing the chain gives
\[
C_{ijk}^{(1,1)}=C_{kji}^{(1,1)}.
\]

\begin{theorem}\label{thm:localization}
Assume \eqref{eq:dimension}. With the notation above,
\begin{equation}\label{eq:localization}
N_2(X)
=
\sum_{0\leq i<j\leq r}C_{ij}^{(2)}
+
\frac12\sum_{j=0}^r
\sum_{\substack{0\leq i,k\leq r\\i\ne j,\ k\ne j}}
C_{ijk}^{(1,1)}.
\end{equation}
The right side is independent of the choice of generic weights $\lambda_0,\ldots,\lambda_r$.
\end{theorem}

\begin{proof}
Apply virtual localization~\cite{GP} to \eqref{eq:lefschetz}. The fixed graphs have either one edge of degree two or two edges of degree one. Proposition~\ref{prop:local} and \eqref{eq:twist} give the corresponding summands.

For a chain with $i\ne k$, the ordered triples $(i,j,k)$ and $(k,j,i)$ represent the same fixed stratum, so the factor $1/2$ removes double counting. When $i=k$, reversal is a label preserving graph automorphism of order two, and the same factor accounts for its stabilizer.

By \eqref{eq:dimension}, the equivariant pushforward has codimension zero and lies in $A_{\mathbb T}^0(\mathrm{point})=\QQ$. It is therefore independent of the generic weights.
\end{proof}

\section{The coefficient formula}\label{sec:hilbert}
The quotient bundle construction \eqref{eq:Em} was used by Katz~\cite{Katz} for the quintic threefold and by Bonavero--H\"oring~\cite{BH} for complete intersections. The interpretation of enumerative counts as top Chern numbers on parameter spaces is discussed in~\cite{FultonIT,EH3264}. By Lemma~\ref{lem:geometry}, the section of 
\[E_X=\bigoplus_{i=1}^s E_{m_i}\] 
has a finite reduced zero scheme parametrizing the embedded conics on $X$. Hence
\begin{equation}\label{eq:conicintegral}
n_2(X)
=
\int_P c_{2m+s}(E_X)
=
\int_P\prod_{i=1}^s c_{2m_i+1}(E_{m_i}).
\end{equation}
Let $x_1,x_2,x_3$ be the Chern roots of $S^\vee$, and write $\mathbf{x}=(x_1,x_2,x_3)$. For each integer $d\geq1$, define
\begin{equation}\label{eq:Cm}
C_d(\mathbf{x};\xi)
=
[t^{2d+1}]
\frac{\displaystyle
\prod_{\substack{a,b,c\geq0\\a+b+c=d}}
\bigl(1+(ax_1+bx_2+cx_3)t\bigr)}
{\displaystyle
\prod_{\substack{a,b,c\geq0\\a+b+c=d-2}}
\bigl(1+(ax_1+bx_2+cx_3-\xi)t\bigr)}.
\end{equation}
The denominator is interpreted as one when $d=1$. By \eqref{eq:Em}, this polynomial represents $c_{2d+1}(E_d)$.

Write
\[
Y_2=\{2x_1,2x_2,2x_3,x_1+x_2,x_1+x_3,x_2+x_3\}
\]
for the multiset of Chern roots of $\Sym^2S^\vee$, and put
\begin{equation}\label{eq:RX}
R_X(\mathbf{x})
=
\sum_{y\in Y_2}
\frac{\prod_{i=1}^s C_{m_i}(\mathbf{x};-y)}
{\prod_{\substack{z\in Y_2\\z\ne y}}(z-y)}.
\end{equation}
This expression is initially interpreted in
$\QQ(x_1,x_2,x_3)$.

Indeed, with our line convention for projectivization, the projective bundle relation is
\[
\prod_{y\in Y_2}(\xi+y)=0.
\]
For a polynomial $F_0(\xi)$ with coefficients symmetric in $x_1,x_2,x_3$, the corresponding pushforward formula is
\[
q_*F_0(\xi)
=
\sum_{y\in Y_2}
\frac{F_0(-y)}
{\prod_{\substack{z\in Y_2\\z\ne y}}(z-y)}.
\]
Equivalently,
\[
q_*(\xi^{5+k})=(-1)^kh_k(Y_2)
\quad(k\geq0),
\]
while $q_*(\xi^j)=0$ for $0\leq j<5$. Here $h_k(Y_2)$ is the complete symmetric polynomial of degree $k$ in $Y_2$. Applying this formula gives
\[
q_*\left(\prod_{i=1}^s C_{m_i}(\mathbf{x};\xi)\right)=R_X(\mathbf{x}).
\]
In particular, the apparent denominators in \eqref{eq:RX} cancel, and $R_X(\mathbf{x})$ is a symmetric homogeneous polynomial of degree
\[
2m+s-5=3(r-2).
\]
Thus \eqref{eq:conicintegral} becomes
\[
n_2(X)=\int_{G_2}R_X(\mathbf{x}).
\]

\begin{theorem}\label{thm:hilbert}
Under the hypotheses of Theorem~\ref{thm:main}, the number of embedded conics on $X$ is
\begin{equation}\label{eq:hilbert}
\begin{aligned}
n_2(X)
&=
\frac1{3!}[x_1^rx_2^rx_3^r]
\left(
R_X(\mathbf{x})\prod_{\substack{1\leq i,j\leq3\\i\ne j}}
(x_i-x_j)
\right)\\
&=
[s_{(r-2,r-2,r-2)}]R_X.
\end{aligned}
\end{equation}
Here $[x_1^rx_2^rx_3^r]$ denotes monomial coefficient extraction, and $[s_{(r-2,r-2,r-2)}]R_X$ denotes the coefficient of the rectangular Schur polynomial $s_{(r-2,r-2,r-2)}(x_1,x_2,x_3)$ in the Schur expansion of $R_X(\mathbf{x})$.
\end{theorem}
\begin{proof}
By \eqref{eq:Em},
\[
c_{2d+1}(E_d)=C_d(\mathbf{x};\xi)
\quad(d\geq1).
\]
Applying the projective bundle pushforward formula to \eqref{eq:conicintegral} gives
\[
n_2(X)
=
\int_{G_2}
q_*\left(\prod_{i=1}^s C_{m_i}(\mathbf{x};\xi)\right)
=
\int_{G_2}R_X(\mathbf{x}).
\]
The polynomial $R_X$ is symmetric and homogeneous of degree
\[
2m+s-5=3(r-2)=\dim G_2.
\]
The Grassmannian coefficient formula~\cite{DP,Hiep2019,HGV} therefore yields
\[
\int_{G_2}R_X(\mathbf{x})
=
\frac1{3!}[x_1^rx_2^rx_3^r]
\left(
R_X(\mathbf{x})\prod_{\substack{1\leq i,j\leq3\\i\ne j}}
(x_i-x_j)
\right).
\]
For the second equality, expand
\[
R_X(\mathbf{x})
=
\sum_{\substack{\lambda\vdash3(r-2)\\\ell(\lambda)\leq3}}
a_\lambda s_\lambda(\mathbf{x}).
\]
The class $s_\lambda(S^\vee)$ vanishes unless $\lambda$ fits inside the $3\times(r-2)$ rectangle. Among partitions of size $3(r-2)$, only $\lambda=(r-2,r-2,r-2)$ fits inside this rectangle. Its Schur class is the class of a point, so
\[
\int_{G_2}s_{(r-2,r-2,r-2)}(S^\vee)=1.
\]
Consequently,
\[
n_2(X)
=
[s_{(r-2,r-2,r-2)}]R_X,
\]
which proves \eqref{eq:hilbert}.
\end{proof}

\begin{proposition}\label{prop:interpolation}
Let $S$ be the tautological subbundle on $\Gr(3,r+1)$, and let $x_1,x_2,x_3$ be the Chern roots of $S^\vee$. For every symmetric homogeneous polynomial $\varphi(\mathbf{x})\in\QQ[x_1,x_2,x_3]$ of degree $3(r-2)$,
\[
\int_{\Gr(3,r+1)}\varphi(S^\vee)
=
\frac16[x_1^rx_2^rx_3^r]
\left(
\varphi(\mathbf{x})\prod_{\substack{1\leq i,j\leq3\\i\ne j}}
(x_i-x_j)
\right).
\]
\end{proposition}

\begin{proof}
Choose pairwise distinct rational numbers $\lambda_0,\ldots,\lambda_r$, and put
\[
Q_0(t)=\prod_{a=0}^r(t-\lambda_a),
\quad
F_0(x)=
\varphi(\mathbf{x})\prod_{\substack{1\leq i,j\leq3\\i\ne j}}
(x_i-x_j).
\]
Since $\deg F_0=3r$, multivariate Lagrange interpolation~\cite[Lemma 1 and Proposition 1]{HGV} gives
\[
[x_1^rx_2^rx_3^r]F_0
=
\sum_{a,b,c=0}^r
\frac{F_0(\lambda_a,\lambda_b,\lambda_c)}
{Q_0'(\lambda_a)Q_0'(\lambda_b)Q_0'(\lambda_c)}.
\]
Terms with repeated indices vanish because of the Vandermonde factors in $F_0$. For a three element subset $I\subset\{0,\ldots,r\}$, write $\varphi(\lambda_I)$ for the evaluation of $\varphi$ at its three weights; symmetry makes this independent of their ordering. For any ordering $(a,b,c)$ of $I$, we have
\[
\prod_{i\in I}Q_0'(\lambda_i)
=
\left(
\prod_{\substack{i,j\in I\\i\ne j}}
(\lambda_i-\lambda_j)
\right)
\left(
\prod_{i\in I}\prod_{j\notin I}
(\lambda_i-\lambda_j)
\right).
\]
The first factor cancels exactly with the Vandermonde product in the numerator. Thus, each ordering contributes
\[
\frac{\varphi(\lambda_I)}
{\prod_{i\in I}\prod_{j\notin I}
(\lambda_i-\lambda_j)}.
\]
There are six orderings of $I$, so
\[
\frac16[x_1^rx_2^rx_3^r]F_0
=
\sum_{\substack{I\subset\{0,\ldots,r\}\\|I|=3}}
\frac{\varphi(\lambda_I)}
{\prod_{i\in I}\prod_{j\notin I}
(\lambda_i-\lambda_j)}.
\]
The right hand side is the Grassmannian localization formula for $\int_{\Gr(3,r+1)}\varphi(S^\vee)$. Indeed, taking $\lambda_0,\ldots,\lambda_r$ as the equivariant Chern roots of the dual ambient vector space, the Chern roots of $S^\vee$ at the fixed point $I$ are $\lambda_i$ for $i\in I$, and the tangent weights are $\lambda_i-\lambda_j$ for $i\in I$, $j\notin I$. This proves the formula.
\end{proof}

\begin{remark}\label{rem:cancel}
Let $\alpha\in\ZZ_{\geq0}^3$ satisfy $|\alpha|=2$, and put $y_\alpha=\alpha\cdot \mathbf{x}$.
For $d\geq2$, substituting $\xi=-y_\alpha$ into \eqref{eq:Cm} changes each denominator factor to
\[
1+((\gamma+\alpha)\cdot \mathbf{x})t,
\quad
\gamma\in\ZZ_{\geq0}^3,\quad |\gamma|=d-2.
\]
The correspondence $\beta=\gamma+\alpha$ identifies these factors with the numerator factors indexed by $\beta\geq\alpha$ componentwise. After cancellation, the remaining product is
\[
\prod_{\substack{\beta\in\ZZ_{\geq0}^3,\ |\beta|=d\\
                 \beta\not\geq\alpha}}
\bigl(1+(\beta\cdot \mathbf{x})t\bigr).
\]
It has exactly
\[
\binom{d+2}{2}-\binom{d}{2}=2d+1
\]
factors. Taking the coefficient of $t^{2d+1}$ therefore gives
\begin{equation}\label{eq:cancel}
C_d(\mathbf{x};-y_\alpha)
=
\prod_{\substack{\beta\in\ZZ_{\geq0}^3,\ |\beta|=d\\
                 \beta\not\geq\alpha}}
(\beta\cdot \mathbf{x}).
\end{equation}
For $d=1$, no vector $\beta$ of weight one satisfies $\beta\geq\alpha$, so the same formula gives 
\[C_1(\mathbf{x};-y_\alpha)=x_1x_2x_3.\] 
Thus every evaluation appearing in \eqref{eq:RX} can be computed directly as a product of linear forms, without expanding a formal power series.
\end{remark}

\begin{corollary}\label{cor:CY}
Let $X$ be a general complete intersection Calabi--Yau threefold in ordinary projective space. Then
\[
N_2(X)=n_2(X)+\frac18n_1(X),
\]
where $n_1(X)$ is the number of lines on $X$, $N_2(X)$ is calculated by \eqref{eq:localization}, and $n_2(X)$ by \eqref{eq:hilbert}. If $n_d^{\mathrm{BPS}}(X)$ denotes the genus zero degree-$d$ BPS invariant, then
\[
n_2^{\mathrm{BPS}}(X)=n_2(X).
\]
\end{corollary}

\begin{proof}
The Calabi--Yau condition gives $m=r+1$. Since $n=r-s=3$, we have $m=s+4$, so \eqref{eq:dimension} holds and $h=0$. By Lemmas~\ref{lem:geometry} and~\ref{lem:normal}, the scheme of lines is finite and reduced, and every line $\ell\subset X$ has the normal bundle
\[
N_{\ell/X}\simeq\OO_\ell(-1)^{\oplus2}.
\]
Every degree one stable map is an isomorphism onto a line. It is unobstructed and has a trivial automorphism group, so
\[
n_1(X)=\#F.
\]
For an embedded smooth conic $C\subset X$, Lemma~\ref{lem:geometry} gives $H^0(C,N_{C/X})=0$. Adjunction gives $\deg N_{C/X}=-2$. Since $N_{C/X}$ has rank two, its splitting on $C\simeq\PP^1$ is therefore
\[
N_{C/X}\simeq\OO_{\PP^1}(-1)^{\oplus2}.
\]
Thus, each embedded conic contributes one to $N_2(X)$. There are no unions of distinct intersecting lines in $X$, by Lemma~\ref{lem:geometry}. The remaining stable maps are the compactified double covers of the isolated lines. For each such line, Proposition~\ref{prop:obstruction} identifies the obstruction bundle with $A^{\oplus2}$, and Lemma~\ref{lem:push} gives
\[
\int_{\rho^{-1}([\ell])}\theta^2=\frac18.
\]
Consequently,
\[
N_2(X)=n_2(X)+\frac18\#F
=n_2(X)+\frac18n_1(X).
\]
Finally, the genus zero multiple cover relation \cite[Chapter 7]{CK} gives
\[
n_1(X)=n_1^{\mathrm{BPS}}(X),
\quad
N_2(X)=n_2^{\mathrm{BPS}}(X)
+\frac18n_1^{\mathrm{BPS}}(X).
\]
Comparison with the preceding formula yields $n_2^{\mathrm{BPS}}(X)=n_2(X)$.
\end{proof}

\begin{remark}\label{rem:pointcheck}
The quotient bundle method also applies to conics satisfying point conditions. As a check, let $X_4\subset\PP^4$ be a general quartic threefold, and fix a general point $p\in X_4$. This example lies outside \eqref{eq:dimension}. Write $\PP^4=\PP(V)$, and let $\ell_p\subset V$ be the one-dimensional subspace corresponding to $p$. Planes containing $p$ are parametrized by
\[
G_p=\Gr(2,V/\ell_p)\simeq\Gr(2,4).
\]
Let $U_p$ be its tautological rank-two subbundle, and let $S_p\subset V\otimes\OO_{G_p}$ be its inverse image, so that
\[
0\longrightarrow\ell_p\otimes\OO_{G_p}
\longrightarrow S_p
\longrightarrow U_p
\longrightarrow0.
\]
Evaluation at $p$ defines a rank-five bundle
\[
K_p
=
\ker\left(
\Sym^2S_p^\vee
\longrightarrow
(\ell_p^\vee)^{\otimes2}\otimes\OO_{G_p}
\right).
\]
Thus
\[
q_p:P_p=\PP(K_p)\longrightarrow G_p
\]
parametrizes plane conics through $p$ and $\dim P_p=4+4=8$. Restrict the bundle $E_4$ to $P_p$, and denote this restriction by $\widehat E_4$. Evaluation at $p$ gives an exact sequence
\[
0\longrightarrow E_4^{(p)}
\longrightarrow\widehat E_4
\longrightarrow
(\ell_p^\vee)^{\otimes4}\otimes\OO_{P_p}
\longrightarrow0,
\]
where $\rank E_4^{(p)}=8$. The quartic equation of $X_4$, which vanishes at $p$, induces a section of $E_4^{(p)}$. The space of quartic equations vanishing at $p$ generates this bundle, so Bertini gives a finite reduced zero scheme for a general pair $(X_4,p)$. Moreover, no line on $X_4$ passes through a general point of $X_4$: the scheme of lines has dimension one, and its universal family has dimension two. Consequently, every conic in this zero scheme is smooth.

Let $x,y$ be the Chern roots of $U_p^\vee$, and put $\xi_p=c_1(\OO_{P_p}(1))$. The roots of $S_p^\vee$ are $0,x,y$, while those of $K_p$ are
\[
Y_p=\{x,y,2x,x+y,2y\}.
\]
The quotient construction gives
\[
c_8(E_4^{(p)})
=
[t^8]\,
\frac{
\displaystyle
\prod_{\substack{k,l\geq0\\1\leq k+l\leq4}}
(1+(kx+ly)t)}
{
\displaystyle
\prod_{\substack{k,l\geq0\\k+l\leq2}}
(1+(kx+ly-\xi_p)t)
}.
\]
Applying the projective bundle pushforward formula, or the cancellation procedure of Remark~\ref{rem:cancel}, yields
\[
\begin{aligned}
(q_p)_*c_8(E_4^{(p)})
={}&216(x^4+y^4)
+1260(x^3y+xy^3)
+2232x^2y^2.
\end{aligned}
\]
The number of conics through $p$ is therefore
\[
\begin{aligned}
\int_{P_p}c_8(E_4^{(p)})
&=
[x^3y^2]\left(
(x-y)(q_p)_*c_8(E_4^{(p)})
\right)\\
&=2232-1260\\
&=972.
\end{aligned}
\]
This agrees with the formula recalled in Bonavero--H\"oring~\cite[Proposition 3.2]{BH}, which gives
\[
\frac{8!}{2^5}-\frac{(4!)^2}{2}=972.
\]
It provides an additional check of the Hilbert scheme calculation with an imposed point condition.
\end{remark}
\section{Examples}\label{sec:examples}
\begin{example}
For the five multidegrees of complete intersection Calabi--Yau threefolds without linear defining equations, formulas \eqref{eq:localization} and \eqref{eq:hilbert} recover the known values below, in agreement with~\cite{Hiep2016, Katz}. Here $n_1(X)$ is the number of lines and $n_2(X)$ is the number of embedded conics.
\begin{center}
\begin{tabular}{lrrr}
\toprule
Multidegree & $n_1(X)$ & $n_2(X)$ & $N_2(X)$\\
\midrule
$(5)$
& $2\,875$ & $609\,250$ & $4\,876\,875/8$\\
$(4,2)$
& $1\,280$ & $92\,288$ & $92\,448$\\
$(3,3)$
& $1\,053$ & $52\,812$ & $423\,549/8$\\
$(3,2,2)$
& $720$ & $22\,428$ & $22\,518$\\
$(2,2,2,2)$
& $512$ & $9\,728$ & $9\,792$\\
\bottomrule
\end{tabular}
\end{center}
Each row satisfies
\[
N_2(X)-n_2(X)=\frac18n_1(X),
\]
as predicted by Corollary~\ref{cor:CY}.

For the quintic threefold $X_5\subset\PP^4$, put $\sigma_i=e_i(x_1,x_2,x_3)$. Evaluating \eqref{eq:RX} using \eqref{eq:cancel} gives
\[
\begin{aligned}
R_{X_5}=50\bigl(
&288\sigma_1^6
+1080\sigma_1^4\sigma_2
+1036\sigma_1^2\sigma_2^2
+144\sigma_2^3\\
&+2042\sigma_1^3\sigma_3
+881\sigma_1\sigma_2\sigma_3
+2366\sigma_3^2
\bigr).
\end{aligned}
\]
The $\sigma_i$ are elementary symmetric polynomials in the Chern roots of $S^\vee$ on $G_2$. The Schur expansion is
\[
\begin{aligned}
R_{X_5}={}&
14\,400s_{(6)}
+126\,000s_{(5,1)}
+343\,400s_{(4,2)}\\
&+513\,900s_{(4,1,1)}
+239\,000s_{(3,3)}\\
&+1\,132\,250s_{(3,2,1)}
+609\,250s_{(2,2,2)}.
\end{aligned}
\]
Since $r=4$, Theorem~\ref{thm:hilbert} extracts the coefficient of $s_{(2,2,2)}$, giving
\[
n_2(X_5)=609\,250.
\]
Together with $n_1(X_5)=2\,875$, this yields
\[
N_2(X_5)
=
609\,250+\frac{2\,875}{8}
=
\frac{4\,876\,875}{8}.
\]
\end{example}

\begin{example}
Let $X=X_8\subset\PP^6$ be a general hypersurface of degree $8$. Then $n=5$, $h=1$, and $\rank B=3$. The dimension condition holds because
\[
2m+s=2\cdot8+1=17=3\cdot6-1.
\]
Corollary~\ref{cor:five} gives
\[
c_1(B)
=
\left(2+\frac{8\cdot7}{2}\right)\sigma
=
30\sigma,
\quad
D_2(X)=\frac{57}{16}\int_F\sigma.
\]
The Fano scheme $F$ is the zero locus of a general section of $\Sym^8U^\vee$ on $\Gr(2,7)$. Its Pl\"ucker degree is
\[
\begin{aligned}
\int_F\sigma
&=
\int_{\Gr(2,7)}c_9(\Sym^8U^\vee)c_1(U^\vee)\\
&=
[x^6y^5]\left(
(x-y)(x+y)\prod_{j=0}^8(jx+(8-j)y)
\right)\\
&=21\,518\,336.
\end{aligned}
\]
In particular, $F$ is nonempty and hence is a smooth projective curve. Therefore
\[
D_2(X)
=
\frac{57}{16}\cdot21\,518\,336
=
76\,659\,072.
\]
Independent evaluations of \eqref{eq:localization} and \eqref{eq:hilbert} give
\[
\begin{aligned}
N_2(X)&=21\,553\,860\,841\,856,\\
n_2(X)&=21\,553\,784\,182\,784.
\end{aligned}
\]
Their difference is
\[
N_2(X)-n_2(X)=76\,659\,072=D_2(X),
\]
verifying the correction formula for a moving one-dimensional family of lines.
\end{example}

\begin{example}
Let $X$ be a general complete intersection of type $(6,3)$ or $(5,4)$ in $\PP^7$. In both cases, $m=9$, $s=2$, and $n=5$ hold, so
\[
2m+s=20=3\cdot7-1,
\quad h=1.
\]
Formula \eqref{eq:fivefactor} gives
\[
D_2(X_{6,3})
=
\frac{1+6\cdot5+3\cdot2}{16}
\int_{F_1(X_{6,3})}\sigma
=
\frac{37}{16}\int_{F_1(X_{6,3})}\sigma,
\]
and
\[
D_2(X_{5,4})
=
\frac{1+5\cdot4+4\cdot3}{16}
\int_{F_1(X_{5,4})}\sigma
=
\frac{33}{16}\int_{F_1(X_{5,4})}\sigma.
\]
The Grassmannian coefficient formula yields the Pl\"ucker degrees
\[
\int_{F_1(X_{6,3})}\sigma=3\,320\,352,
\quad
\int_{F_1(X_{5,4})}\sigma=2\,329\,600.
\]
These positive degrees show that both Fano schemes are nonempty, smooth projective curves. The resulting corrections, together with independent evaluations of \eqref{eq:localization} and \eqref{eq:hilbert}, are
\begin{center}
\begin{tabular}{lrrr}
\toprule
Multidegree & $N_2(X)$ & $n_2(X)$ & $D_2(X)$\\
\midrule
$(6,3)$
& $244\,484\,734\,374$
& $244\,477\,056\,060$
& $7\,678\,314$\\
$(5,4)$
& $108\,786\,674\,000$
& $108\,781\,869\,200$
& $4\,804\,800$\\
\bottomrule
\end{tabular}
\end{center}
In each case, the difference $N_2(X)-n_2(X)$ agrees with the correction calculated from the curve of lines.
\end{example}

\begin{example}
Let $X=X_{11}\subset\PP^8$ be a general hypersurface of degree eleven. Then $n=7$, $h=2$, and $\rank B=4$. The dimension condition holds because
\[
2m+s=2\cdot11+1=23=3\cdot8-1.
\]
Taking the coefficient of $t^2$ in \eqref{eq:main} gives
\begin{equation}\label{eq:seven}
D_2(X)
=
\frac18\int_F
\left(
c_2(B)
-\frac12c_1(B)c_1(\Sym^2U^\vee)
+\frac14h_2(\Sym^2U^\vee)
\right),
\end{equation}
where
\[
h_2(E)=c_1(E)^2-c_2(E)
\]
denotes the complete symmetric class of degree two. In terms of the Chern roots $x,y$ of $U^\vee$,
\[
c_1(\Sym^2U^\vee)=3(x+y),
\quad
h_2(\Sym^2U^\vee)=7x^2+10xy+7y^2.
\]
Evaluating \eqref{eq:ambient}, \eqref{eq:localization}, and \eqref{eq:hilbert} gives
\[
\begin{aligned}
D_2(X)
&=\frac{3\,410\,704\,681\,821\,767}{32},\\
N_2(X)
&=\frac{220\,133\,473\,099\,468\,969\,402\,439}{32},\\
n_2(X)
&=6\,879\,170\,927\,773\,883\,986\,896.
\end{aligned}
\]
These values satisfy
\[
N_2(X)-n_2(X)=D_2(X).
\]
Thus the double cover contribution accounts for the difference between the rational Gromov--Witten invariant and the integer conic count. Its denominator is exactly $32=2^{h+3}$, so the denominator bound in Corollary~\ref{cor:denominators} is attained.
\end{example}

All numerical values above can be computed using exact rational arithmetic. For \eqref{eq:localization}, choose generic rational weights, such as $\lambda_i=3^i$, and evaluate the two graph sums. For \eqref{eq:hilbert}, evaluate $R_X$ at each three element subset of the weights, using \eqref{eq:cancel} to compute its numerator factors, and sum the Grassmannian localization contributions in Proposition~\ref{prop:interpolation}. For \eqref{eq:ambient}, sum over two-element subsets $I=\{i,j\}$, substitute $x=\lambda_i$, $y=\lambda_j$, and extract the coefficient of $t^h$. With our convention that the $\lambda_i$ are weights of the dual ambient vector space, the Chern roots of $Q$ at $I$ are $-\lambda_b$ for $b\notin I$. These calculations independently evaluate $N_2(X)$, $n_2(X)$, and $D_2(X)$. In every example, they satisfy
\[
N_2(X)=n_2(X)+D_2(X).
\]
Their agreement provides numerical checks of the localization factors, the gerbe factor, the twist in the branch construction, and the Chern class expression for the correction. The general identity is proved by the obstruction bundle argument in
Section~\ref{sec:covers}.


\begin{thebibliography}{99}
\bibitem{AM} P.~S. Aspinwall and D.~R. Morrison, Topological field theory and rational curves, \emph{Comm. Math. Phys.} \textbf{151} (1993), no.~2, 245--262, doi:10.1007/BF02096768.
\bibitem{BF} K. Behrend and B. Fantechi, The intrinsic normal cone, \emph{Invent. Math.} \textbf{128} (1997), 45--88.
\bibitem{BH} L. Bonavero and A. H\"oring, Counting conics in complete intersections, \emph{Acta Math. Vietnam.} \textbf{35} (2010), no.~1, 23--30.
\bibitem{CK} D.~A. Cox and S. Katz, \emph{Mirror Symmetry and Algebraic Geometry}, Mathematical Surveys and Monographs, vol.~68, American Mathematical Society, Providence, RI, 1999.
\bibitem{CZ2020} C. Ciliberto and M. Zaidenberg, On Fano schemes of complete intersections, in \emph{Polynomial Rings and Affine Algebraic Geometry}, Springer Proc. Math. Stat., vol.~319, Springer, Cham, 2020, pp.~1--39; arXiv:1903.11294.
\bibitem{CZ2022} C. Ciliberto and M. Zaidenberg, Lines, conics, and all that, \emph{Pure Appl. Math. Q.} \textbf{18} (2022), no.~1, 101--176, doi:10.4310/PAMQ.2022.v18.n1.a4; arXiv:1910.11423.
\bibitem{DP} L. Darondeau and P. Pragacz, Universal Gysin formulas for flag bundles, \emph{Internat. J. Math.} \textbf{28} (2017), 1750077.
\bibitem{EH3264} D.~Eisenbud and J.~Harris, \emph{3264 and All That: A Second Course in Algebraic Geometry}, Cambridge University Press, Cambridge, 2016.
\bibitem{FP} C. Faber and R. Pandharipande, Hodge integrals and Gromov--Witten theory, \emph{Invent. Math.} \textbf{139} (2000), no.~1, 173--199, doi:10.1007/s002229900028.
\bibitem{FantechiP} B. Fantechi and R. Pandharipande, Stable maps and branch divisors, \emph{Compos. Math.} \textbf{130} (2002), no.~3, 345--364, doi:10.1023/A:1014347115536.
\bibitem{FultonIT} W.~Fulton, \emph{Intersection Theory}, second edition, Ergebnisse der Mathematik und ihrer Grenzgebiete, 3.~Folge, vol.~2, Springer-Verlag, Berlin, 1998.
\bibitem{FultonPragacz} W.~Fulton and P.~Pragacz, \emph{Schubert Varieties and Degeneracy Loci}, Lecture Notes in Mathematics, vol.~1689, Springer-Verlag, Berlin, 1998.
\bibitem{GP} T. Graber and R. Pandharipande, Localization of virtual classes, \emph{Invent. Math.} \textbf{135} (1999), 487--518.
\bibitem{Hiep2016} D.~T. Hiep, Rational curves on Calabi--Yau threefolds: verifying mirror symmetry predictions, \emph{J. Symbolic Comput.} \textbf{76} (2016), 65--83, doi:10.1016/j.jsc.2015.12.003.
\bibitem{Hiep2019} D.~T. Hiep, Identities involving (doubly) symmetric polynomials and integrals over Grassmannians, \emph{Fund. Math.} \textbf{246} (2019), 181--191.
\bibitem{HGV} D.~T. Hiep, N.~T.~N. Giao, and N.~T.~M. Van, A characterization of the algebraic degree in semidefinite programming, \emph{Collect. Math.} \textbf{74} (2023), 443--455, doi:10.1007/s13348-022-00358-5.
\bibitem{JNS} M. Jinzenji, I. Nakamura, and Y. Suzuki, Conics on a generic hypersurface, arXiv:math/0412527.
\bibitem{Katz} S. Katz, On the finiteness of rational curves on quintic threefolds, \emph{Compos. Math.} \textbf{60} (1986), 151--162.
\bibitem{Kontsevich} M. Kontsevich, Enumeration of rational curves via torus actions, in \emph{The Moduli Space of Curves}, Progr. Math., vol.~129, Birkh\"auser, 1995, pp.~335--368.
\bibitem{PZ} R. Pandharipande and A. Zinger, Enumerative geometry of Calabi--Yau 5-folds, \emph{Adv. Stud. Pure Math.} \textbf{59} (2010), 239--288.
\bibitem{Voisin} C. Voisin, A mathematical proof of a formula of Aspinwall and Morrison, \emph{Compos. Math.} \textbf{104} (1996), no.~2, 135--151.
\end{thebibliography}
\end{document}